\documentclass[11pt,epsfig]{article}
\usepackage{amsmath}
\usepackage{amsthm}
\usepackage{mathrsfs}
\usepackage{epsfig}
\usepackage{graphicx}
\usepackage{hyperref}
\usepackage{amsfonts}
\usepackage{xcolor}
\usepackage[all,cmtip]{xy}
\allowdisplaybreaks[4]
\usepackage{caption,subcaption}
\usepackage{graphicx}
\usepackage{graphics}
\usepackage{float}

\newtheoremstyle{thmm}{1.5ex plus 1ex minus .2ex}{1.5ex plus 1ex minus .2ex}{\rmfamily}{}{\bfseries}{}{1em}{}
\theoremstyle{thmm}
\newtheorem{theorem}{Theorem}[section]
\newtheorem{lemma}{Lemma}[section]

\newtheorem{remark}{Remark}[section]
\renewenvironment{proof}[1][Proof]{\noindent\textit{#1. } }{\hfill$\square$}
\newcommand{\nn}{\nonumber}
\def \endproof{\vrule height8pt width 5pt depth 0pt}
\def\refe#1{(\ref{#1})}

\def\v{\varepsilon}
\def\e{\epsilon}
\def\x{\bf x}

\def\n{\nu}

\def\p{\phi}

\def\wp{\widetilde p}
\def\u{{\boldsymbol{v}}}
\def\wu{{\widetilde \u}}
\def\wp{{\widetilde p}}

\def\w{{\bf w}}

\def\X{{\bf X}}

\def\x{{\bf x}}
\def\v{{\boldsymbol{v}}}
\def\u{{\boldsymbol{u}}}
\def\e{{\boldsymbol{e}}}

\def\wu{{\widetilde{\u}}}

\def\th{{\boldsymbol{\theta}}}
\def\wth{{\widetilde{\th}}}
\def\wxi{{\widetilde{\xi}}}
\def\n{{\bf n}}
\def\r{{\bf r}}
\def\p{{\bf p}}
\def\q{{\bf q}}

\begin{document}

\title{Optimal analysis of penalized lowest-order mixed FEMs for the Stokes--Darcy model}

\date{}

\author{
\setcounter{footnote}{0}
Luling Cao
\footnote{
School of Mathematics and Data Science, Shaanxi University of Science and Technology, Xi'an, 710021, P. R. China ({\tt lulingcao@163.com}).
The work of the author was partially supported by the Shaanxi Province Natural Science Basic Research Program Youth Project (No.2024JC-YBQN-0058) and Shaanxi Province Postdoctoral Science Foundation.
}
~~ and ~~
Weiwei Sun
\footnote{
Advanced Institute of Natural Science, Beijing Normal University,  Zhuhai, P. R. China, and Guangdong Provincial Key Laboratory of Interdisciplinary Research and Application for Data Science, BNU-HKBU United International College, Zhuhai, 519087, P.R.China ({maweiw@uic.edu.cn}).
The work of the author was partially supported
by the NSFC Key Program (No.12231003) and
Guangdong Provincial Key Laboratory IRADS (2022B1212010006, UIC-R0400001-22).
}
}

\maketitle

\begin{abstract}
\par
This paper is concerned with non-uniform fully-mixed FEMs for dynamic coupled Stokes--Darcy model with the well-known Beavers--Joseph--Saffman (BJS) interface condition.
In particular, a decoupled algorithm with the lowest-order mixed non-uniform FE approximations
(MINI for the Stokes equation and RT0-DG0 for the Darcy equation) and the classical Nitsche-type penalty is studied.
The method with the combined approximation of different orders
is commonly used in practical simulations. However, the optimal error analysis
of methods with non-uniform approximations for the coupled Stokes--Darcy flow model
has remained challenging, although the analysis for uniform approximations has been well done.
The key question is how the lower-order approximation to the Darcy flow influences the accuracy
of the Stokes solution through the interface condition.
In this paper, we prove that the decoupled algorithm provides a truly optimal convergence rate
in $L^2$-norm in spatial direction: $O(h^2)$ for Stokes velocity and $O(h)$ for Darcy flow in the coupled
Stokes--Darcy model. This implies that the lower-order approximation to the Darcy flow does not pollute the accuracy of numerical velocity for Stokes flow.
The analysis presented in this paper is based on a well-designed Stokes--Darcy Ritz projection
and given for a dynamic coupled model.
The optimal error estimate holds for more general combined approximations and more general coupled models, including the corresponding model of steady-state Stokes--Darcy flows
 and the model of coupled dynamic Stokes and steady-state Darcy flows.
Numerical results confirm our theoretical analysis and show that the decoupled algorithm is efficient.
\end{abstract}


\section{Introduction}
\setcounter{equation}{0}
The applications of free flow/porous media flow model are found in various fields of engineering and science, including chemical and petroleum engineering \cite{arbogast2007computational, hanspal2006numerical}, biomedical engineering \cite{d2011robust}, and environmental sciences \cite{discacciati2002mathematical}.
In these contexts, the model is commonly described by
the Stokes (or Navier--Stokes) equations
\begin{equation}
\begin{split}
 \partial_t \u_f - \nabla \cdot \mathbb{T}
& =\boldsymbol{f}_f,
\\
\nabla \cdot \u_f & =0
\end{split}
\label{e-f}
\end{equation}
for the free flow in $\Omega_f \times [0, T]$
and the Darcy equations
\begin{equation}
\begin{split}
\u_p  & = -\boldsymbol{K} \nabla \phi_p,
\\
S_0 \partial_t \phi_p + \nabla \cdot \u_p  & =  f_p
\end{split}
\label{e-p}
\end{equation}
for the porous medium flow in $\Omega_p \times [0, T]$, coupled by the interface conditions
\begin{align}
 \u_f \cdot \n_f &=-\u_p \cdot \n_p,   \label{ConservationMass}
\\
-\n_f \cdot \mathbb{T} \cdot \n_f & = p_f - 2\nu \n_f \cdot \mathbb{D}(\u_f) \cdot \n_f =g_0 \phi_p, \label{BalanceForce}\\
-\n_f \cdot \mathbb{T} \cdot \tau_i = -2 \nu \n_f \cdot \mathbb{D}(\boldsymbol{u}_f) \cdot \tau_i
&= \frac{\alpha \nu \sqrt{d}}{\sqrt{\text{trace}(\Pi)}} \u_f \cdot \tau_i,
~~~~i=1,...,d-1,
\label{BJS}
\end{align}
on $\Gamma: = \overline{\Omega}_f \cap \overline{\Omega}_p$. Here \eqref{BJS} is the well-known Beavers--Joseph--Saffman (BJS) condition, as detailed in \cite{wang1,wang2}.
$\mathbb{T}=2\nu \mathbb{D}(\u_f)-p_f \mathbb{I}$
represents the stress tensor,
$\mathbb{D}(\boldsymbol{u}_f)=\frac{1}{2}(\nabla \boldsymbol{u}_f + \nabla^T \u_f)$ denotes the deformation tensor, $\nu >0$ denotes the kinematic viscosity of the fluid, and
$\boldsymbol{f}_f$ corresponds to a given external force.
The hydraulic conductivity tensor is denoted by $\boldsymbol{K}=\text{diag}(K_1,...,K_d)$, and the eigenvalues of its inverse $\boldsymbol{K}^{-1}$ satisfy
$0 < \lambda_{\min} < \lambda(\boldsymbol{K}^{-1}) < \lambda_{\max} < \infty$.
Moreover, $f_p$ denotes the sink/source term in Darcy region,
$\{\tau_i\}_{i=1}^{d-1}$ denote unit tangential vectors to the interface
$\Gamma$, $d$ is the spatial dimension,
$\alpha >0$ denotes a dimensionless constant depending on the geometrical characteristics of the porous medium, $g_0$ denotes the gravitational acceleration and
$\Pi=g_0^{-1}\boldsymbol{K} \nu$.
$S_0 \geq 0$ denotes the mass storativity coefficient.
For $S_0>0$, the system defined above describes a coupled dynamic Stokes--Darcy model,
while for $S_0=0$, the system reduces to a coupled model of dynamic Stokes and steady-state Darcy flows.

Numerical methods and analysis for the coupled Stokes--Darcy (or Navier--Stokes--Darcy) problem have been extensively investigated \cite{badea2010numerical, cao2014parallel, chen2011parallel, discacciati2007robin, MR2519594, gunzburger2018stokes, lipnikov2014discontinuous, mu2007two, MR4194320, MR4167065} and are typically classified into direct methods for solving the coupled problem and decoupled approaches for solving subproblems independently.
Among these, decoupled methods are popular in finite element (FE) approaches due to their advantages, including computational efficiency, implementation flexibility, and the ability to utilize off-the-shelf efficient solvers for each decoupled problem.
In engineering, to obtain a comprehensively simulation for different physical variables, the fully-mixed methods are widely employed in industrial applications. For Stokes--Darcy equation,  a further unknown (velocity) is added in the Darcy region by the mixed formulation to compute both velocity and pressure in all regions.
Because of the continuity of the Darcy velocity and the discontinuity of the gradient of the pressure in applications, the mixed methods with the Raviart--Thomas FE approximation and DG approximation for Darcy velocity-pressure are more popular in practical simulations.
For fully-mixed methods, the transmission conditions are very important. To deal with the interface conditions, two general numerical approaches are employed:
Lagrange multiplier method and penalized method. The former needs to introduce an extra variable and the latter adds a penalized term in the weak formulation of FE approximation.
A common choice for the penalized method is the Nitsche-type penalty term, $\gamma \langle (\u_{fh} - \u_{ph})\cdot \n_f, (\v_{fh} - \v_{ph})\cdot \n_f \rangle $,
which is often used in computations due to its easy decoupling.
Furthermore, because the governing equations differ for the fluid and porous media regions,
it is natural to employ domain decomposition methods (DDMs). These methods have been extensively studied by many authors, with Robin-type DDMs proving to be particularly efficient for solving steady-state models \cite{chen2011parallel, discacciati2007robin, sun2021domain}.
In \cite{layton2002coupling}, Layton et al. employed the Lagrange multiplier method to impose the interface conditions of the steady-state Stokes--Darcy model and demonstrated optimal error estimates in the energy norm for both the velocity and pressure fields.
Burman and Hansbo \cite{burman2007unified} developed a stabilized method for the steady-state Stokes--Darcy model by incorporating the classical Nitsche penalty, employing uniform (equal-order) fully-mixed FE approximations with the same polynomial order for velocity and pressure in different regions, and provided a comprehensive error analysis in the $H^1$-norm.
Recently Sun et al. \cite{sun2021domain} proposed a Robin-type DDM algorithm
with fully mixed FE approximations.
They proved that the DDM iterative algorithm is convergent and the convergence rate is independent of $h$.
More works on fully-mixed methods can be found in \cite{badia2009unified, camano2015new, discacciati2017conforming, gatica2011analysis, SW, urquiza2008coupling, zhou2021analysis}, most of which are based on uniform approximations (equal-order approximations) or provide error estimates in energy norm for the velocity in both free fluid flow and porous media flow.


The main concern of this paper is non-uniform fully-mixed approximations for the dynamic coupled model,
providing optimal analysis in $L^2$-norm.
In practical applications, different physical variables often require different-order approximations.
Of particular interest is the commonly-used lowest order mixed FE approximations (MINI/RT0-DG0) \cite{brezzi2008mixed, boffi2013mixed} for this coupled model.
However, the lowest-order fully-mixed method is based on a non-uniform approximation
with the approximation property of
$O(h^2)$ for Stokes velocity and $O(h)$ for Darcy velocity/pressure in $L^2$-norm.
It is not sure whether the lower-order approximation to the Darcy flow would pollute the numerical
velocity for the Stokes flow.
Numerical algorithms with non-uniform approximations have been applied extensively
for solving multi-physical problems. The pollution was concerned first in \cite{SW} for an incompressible miscible displacement model, which is governed by a concentration equation
and Darcy equation coupled by the velocity-dependent
diffusion-dispersion tensor and the contraction-dependent viscosity in a single domain.
The underlying Stokes--Darcy model is coupled by the interface conditions, which makes its analysis challenging. Numerical simulations for the Robin-type fully-mixed FEMs were presented  in
\cite{sun2021domain}, which shows that the lower-order approximation to the Darcy flow does pollute the Stokes velocity.
Recently a coupled dynamic model of heat equation/free flow/porous media flow was studied in \cite{zhang2023study},
in which the free flow  is described by Navier--Stokes equation and the porous media flow
by a modified Darcy law
$\partial_t \u_p  + \alpha \u_p = -\nabla \phi_p + f$
(the inclusion of $\partial_t \u_p$ is more artificial, which was debated in literature \cite{mccurdy2019convection, vafai2015handbook}).
A Lagrange multiplier type linearized scheme with the decoupling of flows and heat and the coupling
of Stokes flow and Darcy flow, was proposed, where
the lowest-order fully-mixed FE approximation was used in spatial direction.
They presented optimal error estimates and showed both theoretically and numerically that there is no pollution to the numerical velocity of the Stokes flow.
To the best of our knowledge, no analysis of this lowest-order mixed non-uniform FE approximations has been given for the more realistic coupled Stokes--Darcy flow model \refe{e-f}--\refe{BJS}. Of course,
a decoupled scheme is more interested and more efficient.

In this paper, we focus on a decoupled scheme for the commonly-used Stokes--Darcy model \refe{e-f}--\refe{BJS},
in which the popular lowest-order fully-mixed FE approximations are applied for the discretization in spatial direction and the classical Nitsche-type penalty is used for dealing with the interface condition.
The scheme is efficient since
at each time step, one only needs to solve a Stokes equation and a Darcy equation, respectively.
More important is that for the coupled Stokes--Darcy model, we prove that the scheme provides the optimal error estimates for all physical variables involved, particularly
\begin{align*}
\| \u_f^n - \u_{fh}^n \|_{L^2(\Omega_f)}
\le C (\tau + h^2)
\end{align*}
for Stokes flow, which shows that the lower-order approximation to the Darcy flow does not pollute the numerical velocity of the Stokes flow.
The analysis presented in this paper is based on a well-designed Stokes--Darcy Ritz projection,
which can be viewed as a fully mixed FE approximation to the steady-state model of the coupled Stokes--Darcy flow. Moreover, our theoretical proof is given uniformly for $S_0 \ge 0$. Therefore,
the optimal error estimates presented in this paper
hold also for steady-state Stokes--Darcy model \cite{discacciati2007robin, urquiza2008coupling}
and for coupled dynamic Stokes--steady state Darcy model ($S_0=0$) \cite{LWZ-2020}.

The rest of this paper is organized as follows.
In Section 2, we present a decoupled algorithm with the lowest-order fully-mixed approximation
for the Stokes--Darcy model \refe{e-f}--\refe{BJS}  and our main theoretical results.
In Section 3, we present several useful lemmas and introduce a Stokes--Darcy Ritz
projection, which can be viewed as a steady-state model of the coupled Stokes--Darcy flow.
With the well-designed Ritz projection, we prove the optimal error estimates of the numerical scheme: second-order accuracy for the Stokes velocity and first-order accuracy for others.
In Section 4, we establish the approximation properties of the Stokes--Darcy Ritz projection.
Numerical examples are given in Section 5 to confirm our theoretical results.

\section{Algorithm and main results}
\setcounter{equation}{0}
In this section, we present a decoupled algorithm for the Stokes--Darcy model \refe{e-f}--\refe{BJS}
with the lowest-order mixed non-uniform FE approximations and our main theoretical results.


\begin{figure}[htbp!]
  \centering
  \includegraphics[width=8.25cm]{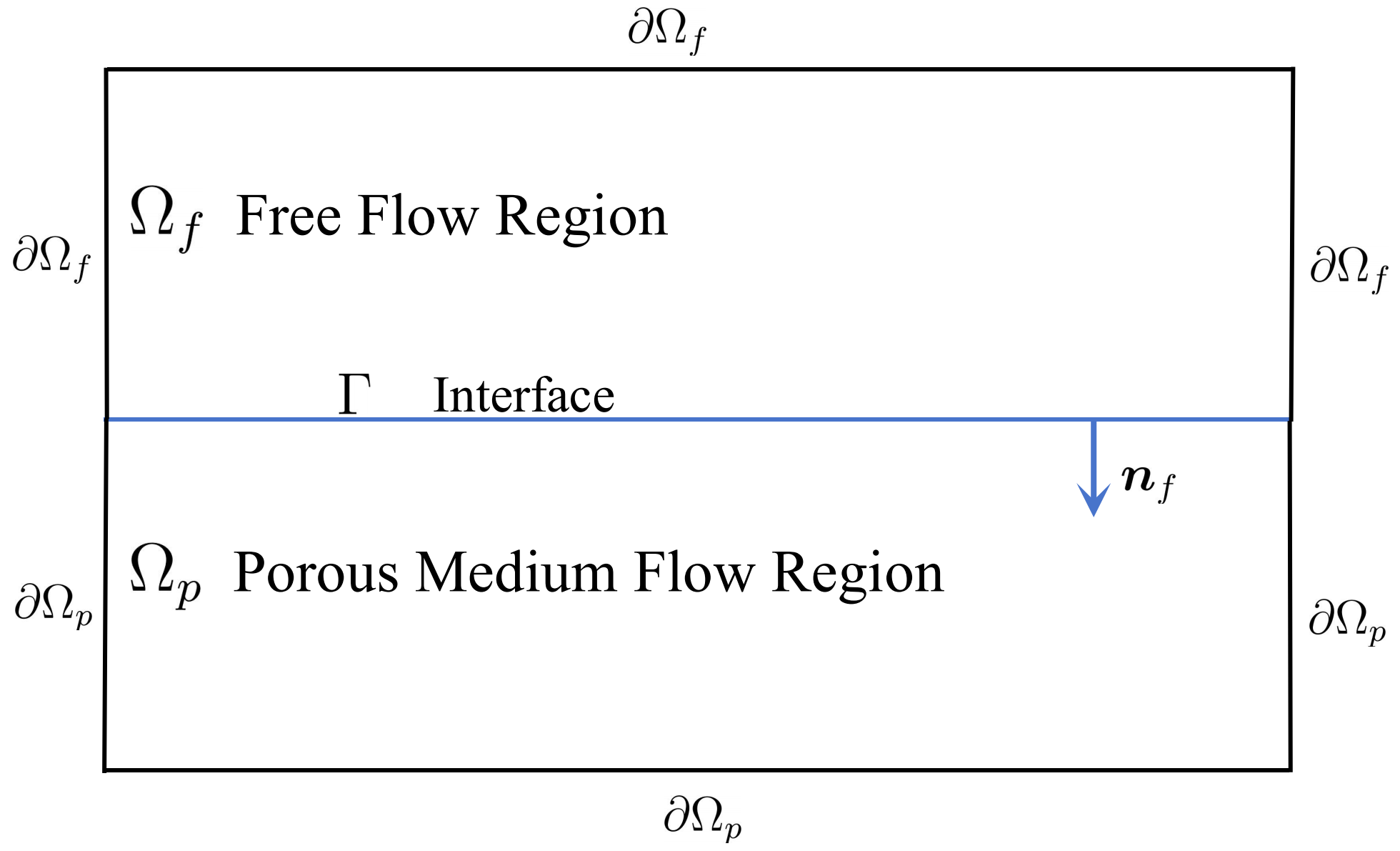}
  \caption{A sketch of the free flow region $\Omega_f$, porous media flow region $\Omega_p$ and the interface $\Gamma$.}\label{region}
\end{figure}

For simplicity, assume that $\Omega \subset \mathbb{R}^d(d=2,3)$ is a bounded domain, divided into two nonintersecting subdomains, $\Omega_f$ and $\Omega_p$, by an interface denoted as $\Gamma$, i.e., $\overline{\Omega}=\overline{\Omega}_f \cup \overline{\Omega}_p$, $\Omega_f \cap \Omega_p=\emptyset$, and $\overline{\Omega}_f \cap \overline{\Omega}_p =\Gamma$ (see Fig.\ref{region}). Let $\n_f$ and $\n_p$ represent the unit outward normal directions on $\partial \Omega_f$ and $\partial \Omega_p$, respectively, with $\n_f=-\n_p$ on $\Gamma$.
The following boundary conditions are imposed:
\begin{align}
  \u_f &=  0 \qquad \quad \, \, \, \, \,  \mbox{ on } \Gamma_f=\partial \Omega_f \backslash \Gamma,
\nonumber\\
 \phi_p & = 0 \qquad \quad  \, \, \, \, \, \mbox{ on } \Gamma_p^D,
\label{Boundary1} \\
  \u_p \cdot \n_p & =0 \qquad \qquad \mbox{ on } \Gamma_p=\partial \Omega_p \backslash \Gamma
\backslash
\Gamma_p^D, \nonumber
\end{align}
and the initial conditions are given by
\begin{align}
\u_f(0,x) = \u_f^0, \qquad \phi_p(0,x) = \phi_p^0
\label{initial}
\end{align}
for Stokes and Darcy equations, respectively.
In addition, we assume that the boundary
$\partial \Omega$ and the interface $\Gamma$ are piecewise smooth.

\subsection{Notation and weak formulation}
Here some notations and functional spaces are introduced.
In what follows, for the subscript $z \in \{f, p\}$, let $L^2(\Omega_z)$ and $H^k(\Omega_z)$ denote the standard Hilbert spaces, with $L^2(\Gamma)$ and $H^k(\Gamma)$
similarly defined on the interface $\Gamma$.
For any functions $u$, $v \in L^{2}(\Omega_z)$ and $w, \xi \in L^2(\Gamma)$, the inner products and norms are defined as
\begin{align}
& (u,\, v)_{\Omega_z}  = \int_{\Omega_z} u(\x) v(\x)\, {\mathrm{d}} \x,
\qquad\, \left \|u\right \|^2_{L^2(\Omega_z)} := (u, \, u)_{\Omega_z} ,
\nn\\
& \langle w,\,  \xi \rangle  = \int_{\Gamma} w(\x) \xi(\x)\, {\mathrm{d}} \x,
\qquad \| w\|_{L^2(\Gamma)}^2 := \langle w, \, w \rangle.
\nn
\end{align}
The simplified notations $(u, v)$, $\| \cdot \|_0$, $\| \cdot \|_1$, $\| \cdot \|_2$ and $\|\cdot\|_{\Gamma}$
will often be used to denote the inner product
and the corresponding $L^2$-norm, $H^1$-norm, $H^2$-norm and $L^2(\Gamma)$-norm, respectively, if not confused.
The functional spaces are defined by
\begin{align*}
& \X_f=\left\{ \v_f \in H^1(\Omega_f)^d: \v_f=0~~\text{on}~\Gamma_f   \right\},
~~~~Q_f= L^2(\Omega_f),
\\
& \X_p=\left\{ \v_p \in H(\mathrm{div},\Omega_p): \v_p \cdot \n_p =0~~\text{on}~\Gamma_p  \right\},
~~~~Q_p=\left \{ \psi \in L^2(\Omega_p):  \psi = 0~~\text{on}~\Gamma_p^D \right \}.
\end{align*}
Also we denote by $\X:=\X_f \times \X_p$ and $Q:=Q_f \times Q_p$ the product spaces for global domain $\Omega$.
The space $\X_p$ is equipped with the norm
$$
\|\v_p\|_{\mathrm{div}}:=\big( \|\v_p\|_{L^2(\Omega_p)}^2 + \|\nabla \cdot \v_p\|_{L^2(\Omega_p)}^2 \big)^{1/2}.
$$
We denote  by $H^{-1}(\Omega_p)$ the dual space of $H^1_{0,D}(\Omega_p):= H^1(\Omega_p)\cap Q_p$,
equipped with the standard dual norm
\begin{align*}
\|v^{\ast}\|_{H^{-1}}=\sup_{0 \neq v \in H^1_{0,D}(\Omega_p)}\frac{( v^{\ast},\, v )}{\|v\|_{1}}.
\end{align*}

The weak formulation of the Stokes--Darcy model \refe{e-f}--\refe{initial} reads as follows: find $\u=(\u_f,\, \u_p) \in \X_f \times \X_p$ with $\u_f \cdot \n_f + \u_p \cdot \n_p =0$ on $\Gamma$ and $\p=(p_f,\,\phi_p) \in Q_f \times  Q_p$
such that
\begin{equation}\label{WeakFormulation}
\begin{split}
( \partial_t \u_f,  \, \v_f )
+a(\u, \, \v)-b(\v, \, \p)
+a_{\Gamma}(\phi_p, \, (\v_f-\v_p)\cdot \n_f)
&=(\boldsymbol{f}_f, \, \v_f)~~\forall\v:=(\v_f, \, \v_p) \in \X,\\
g_0 S_0 (\partial_t \phi_p, \, q_p)
+b(\u, \, \q)&=g_0(f_p, \, q_p)~~\forall\q:=(q_f, \, q_p) \in Q,
\end{split}
\end{equation}
where
\begin{align*}
&a(\u,  \, \v)
=a_f(\u_f, \, \v_f)
+a_p(\u_p, \, \v_p),\\
&a_f(\u_f, \, \v_f)
=2\nu (\mathbb{D}(\u_{f}), \, \mathbb{D}(\v_{f}))
+\frac{\alpha \nu \sqrt{d}}{\sqrt{\mathrm{trace}(\Pi)}} \sum_{i=1}^{d-1} \langle \u_{f} \cdot \tau_i, \, \v_{f} \cdot \tau_i \rangle,
~~
\\
& a_p(\u_p, \, \v_p)
=g_0 \boldsymbol{K}^{-1} (\u_{p}, \, \v_{p}),
\\
&b(\v, \, \p)
=b_f(\v_f, \, p_f) + b_p(\v_p, \, \phi_p),
\\
& b_f(\v_f, \, p_f)=(p_f, \, \nabla \cdot \v_f),
~~~~b_p(\v_p, \, \phi_p)=g_0 (\phi_p, \, \nabla \cdot \v_p),\\
&
a_{\Gamma}(\phi_p, \, (\v_f-\v_p)\cdot \n_f)=g_0\langle \phi_p, \, (\boldsymbol{v}_f -\boldsymbol{v}_p) \cdot \n_f \rangle.
\end{align*}
The well-posedness of the model problem \eqref{WeakFormulation} is well done, $e.g.$, see the literatures \cite{galvis2007non, layton2002coupling}.

To embed the interface condition in \eqref{WeakFormulation}, an equivalent penalized weak formulation is defined by
\begin{align}
(\partial_t \u_f, \, \v_f )
+g_0 S_0 (\partial_t \phi_p,\, q_p)
&+a(\u, \, \v)-b(\v, \, \p)
+b(\u, \, \q)
+a_{\Gamma}(\phi_{p}, \, [\v])
-a_{\Gamma}(q_p, \, [\u])
+\gamma \langle [\u], \, [\v] \rangle
\nn \\
& =(\boldsymbol{f}_f, \, \v_f)
+g_0(f_p, \, q_p), \label{Nitsche-typeFormulation}
\end{align}
where $[\v]: = (\v_f - \v_p) \cdot \n_f$ denotes the jump on the interface.

%
%
\subsection{Fully discrete algorithm with the non-uniform approximations }
Let $T_h=\{K\}$ be a quasi-uniform partition of
$\overline{\Omega}=\overline{\Omega}_f \cup \overline{\Omega}_p$
with the size $h:=\max_{K \in T_h} \{ h_K; h_K=\text{diam}(K)\} >0$.
The element $K \in T_h$ is a triangle for $d=2$ and a tetrahedron for $d=3$.
The subdivisions $T_h^f$ and $T_h^p$ are induced on the regions $\Omega_f$ and $\Omega_p$, respectively, match at the interface $\Gamma$, and form the mesh $T_h := T_h^f \cup T_h^p$.
We use the commonly-used lowest order mixed non-uniform FE approximations (MINI/RT0-DG0) \cite{acosta2011error, MR3584582, brezzi2008mixed, gatica2014simple} and
define the corresponding finite element spaces by
\begin{align*}
&\X_f^h
=\{ \v_{fh} \in \X_f: \v_{fh}|_{K} \in P_{1b}, K \in T_h^f \},~~~~
Q_f^h=
\{ q_{fh} \in Q_f: q_{fh}|_{K} \in P_1 , K \in T_h^f \},\\
&\X_{p}^h=\{ \v_{ph} \in \X_p: \v_{ph}|_{K} \in RT_0 , K \in T_h^p \},~~~~
Q_p^h=\{ q_{ph} \in Q_p: q_{ph}|_{K} \in  P_0 , K \in T_h^p \}
\end{align*}
and the product spaces by
$$
\X^h =\X_f^h \times \X_p^h \subset \X,
\qquad
Q^h = Q_f^h \times Q_p^h \subset Q.
$$
Here
$(\X_f^h, Q_f^h)$ and $(\X_p^h, Q_p^h)$ satisfy the standard inf-sup conditions \cite{brezzi2008mixed, boffi2013mixed}: there exist two constants $\beta_f, \beta_p >0$, independent of $h$, such that
\begin{align*}
\inf_{0\neq q_{fh} \in Q_f^h} \sup_{0 \neq \v_{fh} \in \X_f^h} \frac{b_f(\v_{fh}, \, q_{fh})}{\|\v_{fh}\|_1 \|q_{fh}\|_0} \geq \beta_f,
~~~~\inf_{0\neq q_{ph} \in Q_p^h} \sup_{0 \neq \v_{ph} \in \X_p^h}
\frac{b_p(\v_{ph}, \, q_{ph})}{\|\v_{ph}\|_{\mathrm{div}}\|q_{ph}\|_0} \geq \beta_p.
\end{align*}
In addition, the discrete inf-sup condition is also satisfied by the spaces
$\X_{f,0}^h:=\{\v_f \in \X_f^h: \v_f|_{\partial \Omega_f}=0 \}$ and $Q_{f,0}^h:=\{q_f \in Q_f^h: \int_{\Omega_f} q_f \mathrm{d}\x=0\}$; see \cite{boffi2013mixed, girault2012finite}.

Moreover, we denote by $\boldsymbol{\mathrm{I}}_f$ and $\mathbb{I}_f$ the standard interpolation operators
to the finite element spaces $\X_f^h$ and $Q_f^h$, respectively. By classical interpolation theory,
for all $K \in T_h^f$ we have
\begin{align*}
&\|\v_f - \boldsymbol{\mathrm{I}}_f \v_f\|_{L^2(K)}
+h\|\v_f - \boldsymbol{\mathrm{I}}_f \v_f\|_{H^1(K)} \leq C h^2 \|\v_f\|_{H^2(K)},\\
&\|q_f - \mathbb{I}_f q_f\|_{L^2(K)} \leq Ch \|q_f\|_{H^1(K)}.
\end{align*}
Let  $\boldsymbol{\mathrm{I}}_p: H^1(\Omega_p)^d \rightarrow \X_p^h$ be a projection operator \cite{brezzi2012mixed, layton2002coupling}, defined by
\begin{align}
(\nabla \cdot \boldsymbol{\mathrm{I}}_p \v_{p},\, w_h)=(\nabla \cdot \v_p, \, w_h)~~~~\forall w_h \in Q_p^h.
\label{lpp}
\end{align}
Let $\mathbb{I}_p: Q_p \rightarrow Q_p^h$ be the $L^2$-projection defined by
\begin{align}
(\mathbb{I}_p q_p, w_h)=(q_p, w_h)~~~~\forall w_h \in Q_p^h.
\label{I_p}
\end{align}
For operators $\boldsymbol{\mathrm{I}}_p$ and $\mathbb{I}_p$, by Lemma 4.2 in \cite{MR3647954, layton2002coupling}, we can further obtain
\begin{align*}
&\|\v_p - \boldsymbol{\mathrm{I}}_p \v_{p}\|_{0} \leq C h \|\v_p\|_{1},\\
&\|q_p - \mathbb{I}_p q_p\|_{0} \leq Ch \|q_p\|_{1},
\\
& \|q_p - \mathbb{I}_p q_p\|_{L^2(\Gamma)} \leq Ch \|q_p\|_2.
\end{align*}

Let $\{t_n=n\tau\}_{n=0}^N$ denote a uniform partition of the time interval $[0,T]$ with a step size $\tau=T/N$. We denote $v^n=v(\x, t_n)$ and
\begin{align*}
&d_t v^{n}=\frac{v^{n}-v^{n-1}}{\tau}.
\end{align*}

With the above notations,
a decoupled algorithm is presented below by utilizing the backward Euler scheme in the time direction and the combined lowest-order finite element approximation in the spatial direction.
\vskip0.1in

{\bf Algorithm (Decoupled Scheme)}
Find $(\u_{ph}^{n+1},  \, \phi_{ph}^{n+1}) \in \X_p^h \times  Q_p^h$ such that
for all $\v_{ph} \in \X^h_p, \,  q_{ph} \in Q^h_p$,
\begin{equation}\label{decoupledF}
\begin{split}
&g_0 S_0(d_t \phi_{ph}^{n+1}, \, q_{ph})
+ a_p(\u_{ph}^{n+1}, \, \v_{ph})
- b_p(\v_{ph}, \, \phi_{ph}^{n+1})
+ b_p(\u_{ph}^{n+1},\, q_{ph})
+ a_\Gamma ( \phi_{ph}^{n+1}, \, \v_{ph} \cdot \n_p )
\\
&-\gamma \langle (\u_{fh}^{n} - \u_{ph}^{n+1} )\cdot \n_p, \, \v_{ph} \cdot \n_p \rangle +a_{\Gamma}(q_{ph}, \, (\u_{fh}^n - \u_{ph}^{n+1}) \cdot \n_p)
=g_0(f_p^{n+1}, \, q_{ph})
\end{split}
\end{equation}
and then,
find $(\u_{fh}^{n+1}, \, p_{fh}^{n+1}) \in \X_f^h \times Q_f^h$ such that for all $ \v_{fh} \in \X^h_f, \,  q_{fh} \in Q^h_f$,
\begin{equation}\label{decoupledP}
\begin{split}
(d_t \u_{fh}^{n+1}, \, \v_{fh} )
&+ a_f(\u_{fh}^{n+1}, \, \v_{fh})
- b_f(\v_{fh}, \, p_{fh}^{n+1})
+ b_f(\u_{fh}^{n+1}, \, q_{fh})\\
&+ a_\Gamma (\phi_{ph}^{n+1}, \, \v_{fh} \cdot \n_f)
+ \gamma \langle (\u_{fh}^{n+1}-\u_{ph}^{n+1})\cdot \n_f,  \, \v_{fh} \cdot \n_f \rangle
=(\boldsymbol{f}_f^{n+1}, \, \v_{fh}),
\end{split}
\end{equation}
for $n=0,1,...,N-1$, where $ \u_{fh}^0 = \boldsymbol{\mathrm{I}}_f \u_f^0$ and
$\phi_{ph}^0= \widetilde{\phi}_{ph}^0$.
\begin{remark}
In Algorithm, the penalty parameter $\gamma=\gamma_0$ is chosen as a constant independent of mesh size
$h$ due to the use of non-uniform approximations, while $\gamma=\gamma_0/h$ is commonly adopted for uniform approximations, corresponding to the classical Nitsche method.

In addition, $\widetilde{\phi}_{ph}^0$ denotes the Stokes--Darcy Ritz projection of $\phi_p^0$, which will be defined in \eqref{stokes-darcy-p}.
The decoupled schemes were studied by many authors with either penalized terms or Lagrange multiplier.  Some slightly different schemes can be found in \cite{layton2002coupling, mahbub2023uncoupling}.
\end{remark}

%
%
\subsection{Main results}
This paper focuses on the optimal analysis of the decoupled algorithm
with the commonly-used lowest-order approximations.
We present our main results  in the following theorem.
The proof will be given in the next section.

For simplicity, assume that the system \eqref{e-f}--\eqref{initial} admits a unique solution satisfying
\begin{align}
&\|\u_{f}\|_{L^{\infty}(J;H^2)} + \|p_f\|_{L^2(J;H^1)}
+\|\u_{p}\|_{L^2(J;H^1)} + \|\phi_p\|_{L^{\infty}(J;H^2)}
+\|\partial_t \u_f\|_{L^2(J;H^2)}
+\|\partial_{tt} \u_f\|_{L^2(J;L^2)}
\nn \\
& \quad + \|\partial_t \phi_p\|_{L^{\infty}(J;H^1)}
+\|\partial_{tt} \phi_p\|_{L^2(J;L^2)}
+\|\boldsymbol{f}_f\|_{L^2(J;L^2)}
+\|f_p\|_{L^\infty (J;H^1)}
\leq C_B \, ,
 \label{regularity}
\end{align}
where $C_B$ is a positive constant and $J=[0,T]$.

\begin{theorem}\label{2-1}
Suppose that the system \eqref{e-f}--\eqref{initial} has a unique solution
$(\u_{f}, \u_{p}, p_{f}, \phi_{p})$ satisfying the regularity \eqref{regularity}.
Then
the fully discrete system \eqref{decoupledF}--\eqref{decoupledP} yields a unique solution
$(\u_{fh}^{m}, \u_{ph}^{m}, p^{m}_{fh}, \phi^{m}_{ph}) \in \X^h \times Q^h$, $m=1,...N$. Furthermore, there exists a positive constant $\tau_0$ such that when $\tau\le\tau_0$, the numerical solution satisfies the following error estimates:
\begin{align}
& \max_{1 \le m \le N} \|\u_f(t_m) - \u_{fh}^m\|_{L^2(\Omega_f)}
+ \tau \sum_{n=1}^m\|\u_f(t_n) - \u_{fh}^n \|_{H^1(\Omega_f)}
\leq C_0 (\tau + h^2),
\label{err-f}\\
& \max_{1 \le m \le N} \|\phi_p(t_m) - \phi_{ph}^m\|_{L^2(\Omega_p)}
+ \tau \sum_{n=1}^{m} \|\u_p(t_n) - \u_{ph}^n\|_{L^2(\Omega_p)}
\leq C_0 (\tau + h),
\label{err-p}
\end{align}
where $C_0$ is a positive constant independent of $h$ and $\tau$.
\end{theorem}

\begin{remark}
Clearly the above error estimates are optimal for both Stokes and Darcy velocity/pressure, even though the decoupled algorithm is based on non-uniform approximations of the coupled model.
\end{remark}

Throughout this paper, $C$ and $C_{\epsilon}$ (dependent on $\epsilon$) denote generic constants, and $\epsilon$ represents a small generic constant, independent of $\tau$, $h$, $N$ and $C_0$ in
Theorem \ref{2-1},
which could be different at different occurrences.

%
%
\section{Analysis}
\setcounter{equation}{0}

For given $(\u_f(t), \, p_f(t), \, \u_p(t), \, \phi_p(t)) \in \X_f \times Q_f \times \X_p \times  Q_p, \forall t \in [0,T]$, let $\mathrm{P}_h(\u_f, \, p_f, \,\u_p, \,\phi_p):=(\widetilde{\u}_{fh}, \,\widetilde{p}_{fh}, \,\widetilde{\u}_{ph}, \,\widetilde{\phi}_{ph}) \in \X_f^h \times Q_f^h \times \X_p^h
\times  Q_p^h$ be a Stokes--Darcy Ritz projection defined by
\begin{equation}\label{stokes-darcy-p}
\begin{split}
& a_f(\widetilde{\u}_{fh}-\u_f, \,\boldsymbol{v}_{fh})
 -b_f(\v_{fh}, \,\widetilde{p}_{fh}-p_f)
+a_p(\widetilde{\u}_{ph}-\u_{p}, \, \v_{ph})
-b_p(\v_{ph}, \,\widetilde{\phi}_{ph}-\phi_p)
\\
& \qquad \qquad \qquad \qquad +a_{\Gamma}(\widetilde{\phi}_{ph}-\phi_p, \,[\v_h])
+\gamma \langle [\widetilde{\u}_h-\u], \, [\v_h] \rangle
=0,
 \\
& b_f(\widetilde{\u}_{fh}-\u_{f}, \, q_{fh}) = 0,
\\
& b_p(\widetilde{\u}_{ph}-\u_{p}, \, q_{ph})
-a_{\Gamma}(q_{ph}, \,[\widetilde{\u}_h-\u]) = 0,
\end{split}
\end{equation}
for all $(\v_{fh}, \, q_{fh}, \, \v_{ph}, \, q_{ph}) \in \X_f^h \times Q_f^h \times \X_p^h \times Q_p^h$.
This Stokes--Darcy Ritz projection plays a key role in our analysis.
The corresponding error estimates are presented in the following theorem
and the proof will be given in Section \ref{section4}.

\begin{theorem}\label{T3-1}
For given $(\u_f, \, p_f, \,\u_p, \, \phi_p)  \in \X_f \times Q_f \times \X_p \times  Q_p$, there exists a unique projection
$(\wu_{fh}, \,\widetilde{p}_{fh}, \,\wu_{ph}, \,\widetilde{\phi}_{ph})$ defined in
\eqref{stokes-darcy-p}, which satisfies
\begin{align}
& \| \wu_{fh}-\u_f\|_{L^2(\Omega_f)}
+ h\| \wu_{fh}-\u_f \|_{H^1(\Omega_f)}
\leq C_0^* h^2,
\label{err-sd-1}
\\
& \| \wu_{ph}-\u_p\|_{L^2(\Omega_p)}
+\| \widetilde{\phi}_{ph}-\phi_p\|_{L^2(\Omega_p)}
+ \| \widetilde{\phi}_{ph}-\phi_p\|_{L^2(\Gamma)}  \le C_0^* h,
\label{err-sd-2}
\\
&\| d_t (\wu_{fh}^{n+1}-\u_f^{n+1})\|_{L^2(\Omega_f)}
+  \|\widetilde{\phi}_{ph} - \phi_p \|_{H^{-1}(\Omega_p)}
+ \|d_t (\widetilde{\phi}_{ph}^{n+1} -\phi_p^{n+1} )\|_{H^{-1}(\Omega_p)}
\le C_0^* h^2,
\label{err-sd-3}
\end{align}
where $C_0^*$ is independent of $h$.
\end{theorem}

\begin{remark}
If $(\u_f, \u_p, p_f, \phi_p)$ is a solution of the corresponding steady-state model,
then \refe{stokes-darcy-p} defines its finite element solution in $\X^h \times Q^h$ and
Theorem \ref{3-1} provides the optimal error estimates for the steady-state Stokes--Darcy model.
\end{remark}

Before proving our main theorem, we present some basic estimates in the spaces
$\X_p^h$ and $Q_p^h$.

\begin{lemma}
[\cite{acosta2011error, layton2002coupling}] \label{3-1}
For given $\v_p \in \X_p$,  $\boldsymbol{\mathrm{I}}_p \v_p \in  \X_p^h$ satisfies
\begin{equation}\label{ProDef}
\begin{split}
a_\Gamma(\w_h \cdot \n_p, \,  (\v_p-\boldsymbol{\mathrm{I}}_p \v_p) \cdot \n_p ) = 0 \qquad  \forall  \w_h \in  \X_p^h,
\end{split}
\end{equation}
as well as the continuity bound for $\epsilon >0$,
\begin{align}\label{c-bound}
\| \boldsymbol{\mathrm{I}}_p \v_p \|_{L^2(\Omega_p)} \le C (\|\v_p\|_{H^{\epsilon}(\Omega_p)}
+ \|\nabla \cdot \v_p \|_{L^2(\Omega_p)} ) \, .
\end{align}
\end{lemma}

Let $\chi_{ph} \in Q_p^h $ satisfy the equation
\begin{align} \
b_p(\v_{ph}, \chi_{ph}) - a_\Gamma(\chi_{ph}, \v_{ph} \cdot \n_p) = ( \boldsymbol{F}_1, \v_{ph})  +
( F_2, \nabla \cdot \v_{ph})
+ \langle F_3, \v_{ph} \cdot \n_p \rangle  \quad \forall \v_{ph} \in \X_p^h,
\label{ba-e}
\end{align}
for given $\boldsymbol{F}_1 \in L^2(\Omega_p)^d, F_2 \in L^2(\Omega_p)$, and $F_3 \in L^2(\Gamma)$.
The stability estimate is given below.

We introduce the mixed boundary value problem as follows:
\begin{equation}\label{BVP}
\begin{split}
& \nabla \cdot (\nabla \eta) = \chi_{ph} \qquad ~\mbox{ in } \Omega_p
, \\
& \eta = 0 \qquad \qquad \quad ~~ \,\,\,\,\,\mbox{ on } \Gamma_p^D
, \\
& \nabla \eta \cdot \n_p = 0 \, \,\,\,\, \qquad \quad \mbox{ on } \Gamma_p
, \\
& \nabla \eta \cdot \n_p =  -\chi_{ph} \,\qquad \mbox{ on } \Gamma.
\end{split}
\end{equation}
By (7.28) in Section 7.3 of \cite{lions1970non},
the following regularity bound holds for $\epsilon \leq 1/2$:
\begin{align}
\| \eta \|_{H^{1+\epsilon}}
+ \|\Delta \eta\|_{0}
\le C (\| \chi_{ph} \|_0 + \| \chi_{ph} \|_\Gamma) \, .
\label{reg-2}
\end{align}
Taking $\v_p = \nabla \eta \in \X_p$ and $\v_{ph} = \boldsymbol{\mathrm{I}}_p \v_p$,  by \refe{lpp},
\begin{align}
b_p(\v_{ph}, \,\chi_{ph})=g_0(\nabla \cdot \v_p, \, \chi_{ph})=g_0 \|\chi_{ph}\|_0^2.
\label{bp-1}
\end{align}
On the other hand, by Lemma \ref{3-1},
$$
a_\Gamma (\chi_{ph},\,  \v_{ph} \cdot \n_p )
=a_\Gamma (\chi_{ph},\,  \v_{p} \cdot \n_p)
= -g_0 \| \chi_{ph} \|_\Gamma^2.
$$
It follows that
\begin{align}
\|  \chi_{ph} \|_0^2   +  \|\chi_{ph} \|_\Gamma^2  \le  C(\| {\bf g}_1 \|_0 \| \v_{ph} \|_0 + \| g_2 \|_0
\| \nabla \cdot \v_{ph} \|_0 +
 \| g_3 \|_\Gamma \| \v_{ph} \cdot \n_p\|_\Gamma)
\, .
\nn
\end{align}
By the definition of the projection $\boldsymbol{\mathrm{I}}_p$ (see \refe{lpp} and \refe{ProDef}), the continuity bound \refe{c-bound}, and the auxiliary problem \refe{BVP},
we have
\begin{align*}
& \| \v_{ph} \|_0 =\|\boldsymbol{\mathrm{I}}_p \v_p\|_0 \leq C (\|\v_p\|_{H^{\epsilon}}
+ \|\nabla \cdot \v_p\|_0),
\nn \\
& \| \v_{ph} \cdot \n_p \|_\Gamma  \leq \|\v_p \cdot \n_p\|_\Gamma =  \| \chi_{ph} \|_\Gamma,
\nn \\
& \|\nabla \cdot  \v_{ph} \|_0  \leq \| \nabla \cdot \v_p \|_0 =  \| \chi_{ph} \|_0.
\nn
\end{align*}
With the regularity bound \refe{reg-2}, the following lemma holds.
\begin{lemma}\label{3-2}
Let $\chi_{ph}\in Q_p^h$ satisfy
the equation \refe{ba-e}.
Then it holds that
\begin{align}
\| \chi_{ph} \|_0  + \|\chi_{ph} \|_\Gamma \le C (\| \boldsymbol{F}_1 \|_0  +  \| F_2 \|_0 + \| F_3 \|_\Gamma )\, .
\label{ba-est-1}
\end{align}
\end{lemma}

\begin{lemma}[\cite{sun2022new}]\label{3-3}
Let $w_h \in Q_p^h$ satisfy the equation
\begin{align*}
(w_h,\, \nabla \cdot \v_h)=(\boldsymbol{f}, \, \v_h)~~~~\forall \v_h \in \X_p^h \cap \{ \v_h \cdot \n_p|_{\Gamma}=0 \},
\end{align*}
for a given $\boldsymbol{f} \in L^2(\Omega_p)^d$.
 Then for any $g_p \in L^2(\Omega_p)$,
\begin{align}\label{H-1L2}
(g_p,\, w_h) \leq C \|\boldsymbol{f} \|_0 (\|g_p\|_{H^{-1}}
+ h \|g_p\|_0 ).
\end{align}
\end{lemma}

Now we turn back to the proof of Theorem \ref{2-1}.
In terms of the above Stokes--Darcy Ritz projection \eqref{stokes-darcy-p},
the error functions are split into two parts:
\begin{align*}
&\u_{fh}^n-\u_{f}^n
=(\u_{fh}^n-  \wu_{fh}^n)
+( \wu_{fh}^n-\u_f^n)
=\th_{fh}^n+\wth_f^n \, , \\
&p_{fh}^n-p_f^n
=(p_{fh}^n- \widetilde p_{fh}^n )+( \widetilde p_{fh}^n - p_f^n)
=\xi_{fh}^n + \wxi_f^n\, ,\\
&\u_{ph}^n - \u_p^n
=(\u_{ph}^n -  \wu_{ph}^n)
+( \wu_{ph}^n - \u_p^n)
=\th_{ph}^n + \wth_p^n \, ,\\
&\phi_{ph}^n - \phi_p^n
=(\phi_{ph}^n - \widetilde \phi_{ph}^n) + (\widetilde \phi_{ph}^n - \phi_p^n)
=\xi^n_{ph} + \wxi_p^n \, .
\end{align*}
Then we rewrite the weak formulation \eqref{WeakFormulation} into
\begin{equation}\label{p-WeakFormulation}
\begin{split}
& ( d_t \u_f^{n+1}, \, \v_f )
+
a(\u^{n+1}, \, \v)-b(\v, \, \p^{n+1})
+a_{\Gamma}(\phi_p^{n+1}, \, (\v_f - \v_p) \cdot \n_f)
\\
& ~~~~~~~~~~~~~~~~~~~~~~~~~~~~~~~ =(\boldsymbol{f}_f^{n+1},\, \v_f)
+ (T_f^{n+1}, \, \v_f)
~~~~~~~\forall~\v:=(\v_f, \, \v_p) \in \X,
\\
& g_0 S_0( d_t \phi_p^{n+1}, \, q_p )
+b(\u^{n+1}, \, \q)=g_0(f_p^{n+1}, \, q_p)
 + g_0 S_0 (T_p^{n+1}, \, q_p) \qquad
 \forall~\q:=(q_f, \, q_p) \in Q,
\end{split}
\end{equation}
where $T_f^n = d_t \u_f^n - \partial_t \u_f^n$ and $T_p^n= d_t \phi_p^n - \partial_t \phi_p^n$ denote the truncation errors, satisfying
\begin{align}
\sum_{n=1}^N \tau  \| T_f^n \|_0^2 \le C \tau^2 \int_0^T \|\partial_{tt} \u_f \|_0^2 \mathrm{d} t  \, ,
~~~~~~\sum_{n=1}^N \tau \| T_p^n \|_0^2 \le C \tau^2 \int_0^T \|\partial_{tt} \phi_p
\|_0^2 \mathrm{d}t \, .
\label{t-bound}
\end{align}
In terms of the above splitting and the Stokes--Darcy Ritz projection, we subtract the two equations
in \refe{p-WeakFormulation} from \eqref{decoupledF} and \eqref{decoupledP}, respectively and add
them together to get the error equation
\begin{align}
&(d_t \th_{fh}^{n+1}, \, \v_{fh}) +a_f(\th_{fh}^{n+1}, \, \v_{fh})
-b_f(\v_{fh}, \, \xi_{fh}^{n+1})  +b_f(\th_{fh}^{n+1}, \, q_{fh})
\nonumber\\
&+g_0 S_0(d_t \xi_{ph}^{n+1}, \, q_{ph})
+a_p(\th_{ph}^{n+1}, \, \v_{ph})
-b_p(\v_{ph}, \, \xi_{ph}^{n+1})
+b_p(\th_{ph}^{n+1}, \, q_{ph}) \nonumber\\
&+a_{\Gamma}(\xi_{ph}^{n+1}, \, [\v_{h}])
-a_{\Gamma}(q_{ph}, \, (\u_{fh}^n - \wu_{fh}^{n+1} - \th_{ph}^{n+1}) \cdot \n_f ) \nonumber\\
&-\gamma \langle ( \u_{fh}^{n}-\u_{fh}^{n+1})\cdot \n_f,  \, \v_{ph} \cdot \n_f \rangle
 +\gamma \langle ( \th_{fh}^{n+1}-\th_{ph}^{n+1})\cdot \n_f, \, [\v_{h}] \rangle
\nonumber\\
&=
-\left(  T_f^{n+1}, \, \v_{fh}  \right)
-g_0 S_0 \left(  T_p^{n+1}, \, q_{ph}  \right)
-\left( d_t \widetilde{\boldsymbol{\theta}}_f^{n+1}, \, \v_{fh} \right)
-g_0 S_0 \left(  d_t \widetilde{\xi}_{p}^{n+1}, \, q_{ph}  \right). \label{Error1}
\end{align}

For $n\geq 0$, take  $(\v_{fh}, \,\v_{ph}, \,q_{fh}, \,q_{ph})=2\tau (\th_{fh}^{n+1}, \,\th_{ph}^{n+1}, \,\xi_{fh}^{n+1}, \,\xi_{ph}^{n+1})$
in \eqref{Error1} to obtain
\begin{align}
& \|\th_{fh}^{n+1}\|_0^2   - \|\th_{fh}^n\|_0^2  + \|\th_{fh}^{n+1}-\th_{fh}^n\|_0^2
+ g_0 S_0 (\|\xi_{ph}^{n+1}\|_0^2  - \|\xi_{ph}^n\|_0^2  +
\|\xi_{ph}^{n+1}-\xi_{ph}^n\|_0^2)
\nn \\
&
\quad +2\tau a_f(\th_{fh}^{n+1}, \, \th_{fh}^{n+1})
+2\tau a_p(\th_{ph}^{n+1}, \, \th_{ph}^{n+1})
+ 2 \gamma \tau
\| (\th_{fh}^{n+1}-\th_{ph}^{n+1}) \cdot \n_f \|_\Gamma^2
\nn \\
& \le \epsilon \tau (\| \th_{fh}^{n+1} \|_0^2
+ S_0\|\xi_{ph}^{n+1}\|_0^2)
+ C_\epsilon \tau (\| T_f^{n+1} \|_0^2
+ S_0 \| T_p^{n+1} \|_0^2)
\nn \\
&\quad
-2 \tau \left( d_t \widetilde{\boldsymbol{\theta}}_f^{n+1}, \, \boldsymbol{\theta}_{fh}^{n+1} \right)
-2 g_0 S_0 \tau \left(  d_t \widetilde{\xi}_{p}^{n+1}, \, \xi_{ph}^{n+1}  \right)
\nn
\\
&\quad
- 2\gamma \tau \langle  (\u_{fh}^{n+1}-\u_{fh}^{n}) \cdot \n_f,  \, \th_{ph}^{n+1} \cdot \n_f \rangle
 - 2\tau
a_{\Gamma}(\xi_{ph}^{n+1}, \, (\u_{fh}^{n+1}-\u_{fh}^n) \cdot \n_f )
\nn
\\
&:= \epsilon \tau \| \th_{fh}^{n+1} \|_0^2
+ C S_0 \tau \|\xi_{ph}^{n+1}\|_0^2
+ C_\epsilon \tau \| T_f^{n+1} \|_0^2
+ C_{\epsilon} S_0 \tau \|T_p^{n+1}\|_0^2
+ \sum_{i=1}^4 J_i^{n+1} \, .
\label{err-1}
\end{align}
By Cauchy--Schwarz inequality and Theorem \ref{T3-1}, we have
\begin{align}
|J_1^{n+1}|=\big| 2 \tau \left( d_t \widetilde{\boldsymbol{\theta}}_f^{n+1}, \, \boldsymbol{\theta}_{fh}^{n+1} \right) \big|
\leq C_{\epsilon} \tau \|d_t \widetilde{\boldsymbol{\theta}}_f^{n+1}\|_0^2 +
\epsilon \tau \|\boldsymbol{\theta}_{fh}^{n+1}\|_0^2
\le
\epsilon \tau \|\boldsymbol{\theta}_{fh}^{n+1}\|_0^2  + C_\epsilon \tau h^4.
\label{J1}
\end{align}
To estimate $J_2^{n+1}$,
letting $\v_{fh}=0, q_{fh}=0, q_{ph}=0, \v_{ph} \cdot \n_f|_{\Gamma}=0$, it follows from \eqref{Error1} that
\begin{align*}
b_p(\v_{ph}, \, \xi_{ph}^{n+1})
=a_p(\th_{ph}^{n+1}, \, \v_{ph}).
\end{align*}
Using Lemma \ref{3-3},
\begin{align*}
|J_2^{n+1}|=\big| 2 g_0 S_0 \tau  \left(  d_t \widetilde{\xi}_{p}^{n+1}, \, \xi_{ph}^{n+1}  \right) \big|
&\leq C S_0 \tau  \|\th_{ph}^{n+1}\|_0 ( \|d_t \widetilde{\xi}_{p}^{n+1}\|_{H^{-1}}
+  h \|d_t \widetilde{\xi}_{p}^{n+1}\|_0 ),
\end{align*}
and by Theorem \ref{T3-1}, this further implies
\begin{align*}
|J_2^{n+1}| \leq \epsilon S_0 \tau \|\th_{ph}^{n+1}\|_0^2 + C_\epsilon S_0 \tau h^4 \, .
\end{align*}
By noting the fact that $\| d_t \wth_f^{n+1} \|_\Gamma  \le \| d_t \wth_f^{n+1} \|_0^{1/2}
\| d_t \wth_f^{n+1} \|_1^{1/2} \le  Ch^{3/2}$ and $\| d_t \u_f^{n+1} \|_\Gamma \le C$, we have
\begin{align*}
\|(\u_{fh}^{n+1}-\u_{fh}^{n})\cdot \n_f\|_{\Gamma}
&\leq C (\|  (\th_{fh}^{n+1} - \th_{fh}^{n} )\cdot \n_f \|_\Gamma
+ \|  (\wth_{f}^{n+1} - \wth_{f}^{n} )\cdot \n_f \|_\Gamma
+ \|  (\u_{f}^{n+1} - \u_{f}^{n} )\cdot \n_f \|_\Gamma)
\nn \\
& \leq C (\|  \th_{fh}^{n+1} \cdot \n_f \|_\Gamma  +\| \th_{fh}^{n} \cdot \n_f \|_\Gamma
+\tau  h^{3/2}
+\tau ),
\end{align*}
which with Theorem \ref{3-1}, we further shows
\begin{align}
| J_3^{n+1}  |
& \le
C \tau \big( \|  \th_{fh}^{n+1} \cdot \n_f \|_\Gamma + \|  \th_{fh}^{n} \cdot \n_f \|_\Gamma
+\tau  h^{3/2}
+\tau \big) (\| (\th_{fh}^{n+1}-\th_{ph}^{n+1})\cdot \n_f \|_\Gamma
+\|\th_{fh}^{n+1} \cdot \n_f \|_\Gamma)\, .
\label{J3}
\end{align}
Similarly,
\begin{align}
|J_4^{n+1}|
\le C \tau \| \xi_{ph}^{n+1} \|_\Gamma \big( \|  \th_{fh}^{n+1} \cdot \n_f \|_\Gamma + \|  \th_{fh}^{n} \cdot \n_f \|_\Gamma
+\tau  h^{3/2}
+\tau  \big).
\label{J4}
\end{align}
On the other hand, taking $\v_{fh}=0, q_{fh}=0, q_{ph}=0$
in \refe{Error1}, we see that
\begin{equation}
\begin{split}
b_p(\v_{ph}, \, \xi_{ph}^{n+1})
- a_{\Gamma} (\xi_{ph}^{n+1}, \, \v_{ph} \cdot \n_p) &= a_p(\th_{ph}^{n+1}, \, \v_{ph})
+ \gamma \langle (\u_{fh}^n-\u_{fh}^{n+1})\cdot \n_f, \, \v_{ph} \cdot \n_p \rangle\\
&~~~~ + \gamma \langle (\th_{fh}^{n+1}-\th_{ph}^{n+1})\cdot \n_f, \, \v_{ph} \cdot \n_p \rangle.
\end{split}
\nn
\end{equation}
Taking  $\chi_{ph}=\xi_{ph}^{n+1}$
in \eqref{ba-e} and applying Lemma \ref{3-2} for the last equation gives
\begin{align}
\| \xi_{ph}^{n+1} \|_0
+ \| \xi_{ph}^{n+1} \|_\Gamma
\le C (\| \th_{ph}^{n+1} \|_0
+ \|(\u_{fh}^{n+1}-\u_{fh}^{n})\cdot \n_f\|_{\Gamma}
+ \|(\th_{fh}^{n+1}-\th_{ph}^{n+1})\cdot \n_f\|_{\Gamma}).
\nn
\end{align}
It follows that
\begin{align}
| J_3^{n+1}| + |J_4^{n+1}|
& \le
C \tau \Big( \|  \th_{fh}^{n+1} \cdot \n_f \|_\Gamma + \| \th_{fh}^{n} \cdot \n_f \|_\Gamma
+\tau  h^{3/2}
+\tau   \Big)
\nn \\
&~~~~\cdot \Big(\| \th_{ph}^{n+1} \|_0
+\| (\th_{fh}^{n+1}-\th_{ph}^{n+1})\cdot \n_f \|_\Gamma
+\|\th_{fh}^{n+1} \cdot \n_f \|_\Gamma
+ \|  \th_{fh}^{n} \cdot \n_f \|_\Gamma  + \tau h^{3/2} + \tau \Big)
\nn \\
&\leq C \tau \Big( \|\th_{fh}^{n+1}\|_0^{1/2}\|\th_{fh}^{n+1}\|_1^{1/2}
+\|\th_{fh}^{n}\|_0^{1/2}\|\th_{fh}^{n}\|_1^{1/2}
+\tau h^{3/2}
+\tau \Big)
\nn \\
&~~~~\cdot \Big(\| \th_{ph}^{n+1} \|_0
+\| (\th_{fh}^{n+1}-\th_{ph}^{n+1})\cdot \n_f \|_\Gamma
+\|\th_{fh}^{n+1}\|_0^{1/2} \|\th_{fh}^{n+1}\|_1^{1/2} \nn \\
&~~~~~~~~+  \|\th_{fh}^{n}\|_0^{1/2} \|\th_{fh}^{n}\|_1^{1/2}
+\tau h^{3/2}
+\tau \Big)
\nn \\
&\leq \epsilon \tau ( \| \th_{ph}^{n+1}  \|_0^2
+  \| \th_{fh}^{n+1} \|_1^2 +   \| \th_{fh}^{n}\|_1^2
+  \|(\th_{fh}^{n+1}-\th_{ph}^{n+1})\cdot \n_f\|_{\Gamma}^2)
\nn \\
&~~~~+ C_\epsilon \tau ( \| \th_{fh}^{n+1}  \|_0^2
+   \| \th_{fh}^{n} \|_0^2)
+ C \tau^3,
\label{J34}
\end{align}
where the trace theorem has been used.

With the above estimates and choosing an appropriate $\epsilon$, \refe{err-1} reduces to
\begin{align*}
 \|\th_{fh}^{n+1}\|_0^2   - \|\th_{fh}^n\|_0^2  & + \|\th_{fh}^{n+1}-\th_{fh}^n\|_0^2
+ g_0 S_0 (\|\xi_{ph}^{n+1}\|_0^2 - \|\xi_{ph}^n\|_0^2 +
\|\xi_{ph}^{n+1} - \xi_{ph}^n\|_0^2 )
\nn \\
&
\quad
+\tau \|\th_{fh}^{n+1} \|_1^2
+\tau \| \th_{ph}^{n+1} \|_{0}^2 + \tau
\| (\th_{fh}^{n+1}-\th_{ph}^{n+1}) \cdot \n_f \|_\Gamma^2
\nn \\
 \le & C_\epsilon \tau (\| \th_{fh}^{n+1} \|_0^2 + \| \th_{fh}^{n} \|_0^2
 + S_0 \|\xi_{ph}^{n+1}\|_0^2)
+ \epsilon \tau \|\th_{fh}^{n}\|_1^2
\nn \\
& + C \tau \| T_f^{n+1} \|_0^2
+ C_\epsilon S_0 \tau \| T_p^{n+1} \|_0^2
+ C \tau^3
+ C \tau h^4,  \qquad \mbox{ for } n\geq 0 \, .
\end{align*}
Summing up the last inequality from $n=0$ to $m-1$ and applying Gronwall inequality, when $\tau < \tau_0$ for some $\tau_0 >0$, the following estimate holds
\begin{align}
& \|\th_{fh}^{m}\|_0^2
+ S_0 \|\xi_{ph}^{m}\|_0^2
+\tau \sum_{n=0}^{m-1} \|\th_{fh}^{n+1} \|_1^2
+\tau \sum_{n=0}^{m-1} \| \th_{ph}^{n+1} \|_{0}^2
\nn \\
&\le \|\th_{fh}^{0}\|_0^2
+ S_0 \|\xi_{ph}^{0}\|_0^2
+  \epsilon \tau \|\th_{fh}^{0}\|_1^2
+C \tau \sum_{n=0}^{m-1} (\| T_f^{n +1} \|_0^2 + S_0 \| T_p^{n + 1}\|_0^2)
+C (\tau^2 + h^4)
\nn \\
&\le C_0(\tau^2 + h^4),
\label{err-3}
\end{align}
with $C_0=\exp\left(C_{\epsilon}T / (1-C_{\epsilon} \tau) \right)$,
where we have noted $\|\th_{fh}^{0}\|_l = \| \boldsymbol{\mathrm{I}}_f \u_f^0 - \wu_{fh}^0 \|_l \le C h^{2-l} (l=0,1)$, $\xi_{ph}^{0}=0$ and \refe{t-bound}.
Combining the last equation and the error estimates presented in Theorem \ref{3-1} completes the proof of Theorem \ref{2-1}.
\endproof

%
%
\section{Proof of Theorem \ref{3-1}}\label{section4}
\setcounter{equation}{0}
To prove Theorem \ref{T3-1}, we introduce the following splitting of error functions:
\begin{align*}
&\wu_{fh}-\u_{f}
=(\wu_{fh}- \boldsymbol{\mathrm{I}}_f \u_{f})
+( \boldsymbol{\mathrm{I}}_f \u_{f}-\u_f )
=\e_{fh}+\e_f,\\
&\wp_{fh}-p_f
=(\wp_{fh}-\mathbb{I}_f p_f )+(\mathbb{I}_f p_f - p_f)
=\eta_{fh} + \eta_f,\\
&\wu_{ph} - \u_p
=(\wu_{ph} - \boldsymbol{\mathrm{I}}_p \u_p)
+(\boldsymbol{\mathrm{I}}_p \u_p - \u_p)
=\e_{ph} + \e_p,\\
&\widetilde \phi_{ph} - \phi_p
=(\widetilde \phi_{ph} - \mathbb{I}_p \phi_p) + (\mathbb{I}_p \phi_p - \phi_p)
=\eta_{ph} + \eta_p.
\end{align*}
Denote
\begin{align}
\e_h =( \e_{fh}, \, \e_{ph}), \quad \e =( \e_f, \, \e_p),  \quad \eta_h=(\eta_{fh}, \, \eta_{ph}), \quad
\eta = (\eta_f, \, \eta_p),
\nn
\end{align}
and
$$
\wu_h=(\wu_{fh}, \,\wu_{ph}) \in \X^h,  \qquad
{\bf \wp}_h=(\widetilde p_{fh}, \, \widetilde \phi_{ph}) \in Q^h \, .
$$
\refe{stokes-darcy-p} can be rewritten into
\begin{equation}\label{Orthogonality}
\begin{split}
a(\wu_h-\u, \,\v_h)
-b(\v_h, \, {\bf \wp}_{h}-\p)
&+b(\wu_h-\u, \,\q_h)
+a_{\Gamma}(\widetilde \phi_{ph}-\phi_{p}, \, [\v_h])\\
&-a_{\Gamma}(q_{ph}, \, [\wu_h-\u])
+\gamma \langle [\wu_h-\u], \, [\v_h] \rangle
=0,
\end{split}
\end{equation}
for any $\v_h=(\v_{fh}, \,\v_{ph}) \in \X^h$ and $\q_h=(q_{fh}, \,q_{ph}) \in Q^h$.

By noting the definition of the splitting errors,
\eqref{Orthogonality} reduces to
\begin{align}
&a(\e_{h}, \,\v_h)
-b(\v_h, \, \eta_h)
+b(\e_h, \, \q_h)
+a_{\Gamma}(\eta_{ph},  \, [\v_h])
-a_{\Gamma}(q_{ph}, \,  [\e_h])
+\gamma \langle [\e_h], \,  [\v_h] \rangle
\nn \\
&=-\big( a(\e, \, \v_h)
-b(\v_h, \, \eta)
+b(\e, \, \q_h)
+a_{\Gamma}(\eta_{p}, \, [\v_h])
-a_{\Gamma}(q_{ph}, \, [\e])
+\gamma \langle [\e], \, [\v_h] \rangle  \big),
\label{err-0}
\end{align}
for all $\v_h=(\v_{fh}, \,\v_{ph}) \in \X^h, \q_h=(q_{fh}, \,q_{ph}) \in Q^h$.
We take $\v_h=\e_h$ and $\q_h=\eta_h$ in \eqref{err-0} to get
\begin{align}
& a(\e_h, \,\e_h)
+\gamma \langle [\e_h], \, [\e_h] \rangle
\nn \\
&=-\big( a(\e, \,\e_h)
-b(\e_h, \,\eta)
+b(\e, \,\eta_h)
+a_{\Gamma}(\eta_{p}, \, [\e_h])
-a_{\Gamma}(\eta_{ph}, \, [\e])
+\gamma \langle [\e], \,[\e_h] \rangle \big)
:=\mathrm{RHS}.
\label{error-1}
\end{align}
It is easy to see
\begin{align}
a(\e_h, \,\e_h)
+\gamma \langle [\e_h], \, [\e_h] \rangle
\geq C(\|\e_{fh} \|_1^2 + \|\e_{ph}\|_{0}^2
+ \|[\e_h]\|_{\Gamma}^2).
\label{a+gamma}
\end{align}
Then we need to estimate RHS.
By noting approximation properties of interpolation/projection operators,
RHS can be bounded by
\begin{align}\label{RHS}
\mathrm{RHS} &\leq Ch(\|\e_{fh}\|_1
+\|\e_{ph}\|_{0}
+ \|\eta_{fh}\|_0
+ \|[\e_h]\|_{\Gamma}) \nn\\
&~~~~+|b_p(\e_{ph},\eta_p)|
+|b_p(\e_p,\, \eta_{ph})|
+|a_{\Gamma}(\eta_{ph}, \, [\e])|
+\gamma |\langle [\e], \, [\e_h] \rangle|.
\end{align}

Since $\mathbb{I}_p \phi_p$ and $\boldsymbol{\mathrm{I}}_p \u_p$ denote the $L^2$-projections of $\phi_p$ and $\u_p$, respectively, we obtain
\begin{equation}\label{bp}
\begin{split}
& b_p(\e_{ph}, \, \eta_p) = (\nabla \cdot \e_{ph}, \, \mathbb{I}_p \phi_p - \phi_p) = 0,
\\
& b_p(\e_p,\, \eta_{ph}) = (\nabla \cdot (\boldsymbol{\mathrm{I}}_p \u_p-\u_p), \, \eta_{ph}) = 0,
\end{split}
\end{equation}
where it is noted that $ \nabla \cdot \e_{ph}, \eta_{ph} \in Q_p^h$.
Then, we estimate the last two terms in \eqref{RHS}.
By noting $ (\u_f - \u_p) \cdot \n_f = 0$ on $\Gamma$, and
\begin{align}\label{IpProInt}
\int_{\Gamma_j} (\u_p \cdot \n_f - (\boldsymbol{\mathrm{I}}_p \u_p) \cdot \n_f) \mathrm{d}s = 0.
\end{align}
At each edge/face $\Gamma_j$ of the interface $\Gamma$, we have
\begin{align}
(\boldsymbol{\mathrm{I}}_p \u_p) \cdot \n_f  = \frac{1}{|\Gamma_j|} \int_{\Gamma_j} \u_p \cdot \n_f \mathrm{d}s =
\frac{1}{|\Gamma_j|} \int_{\Gamma_j} \u_f \cdot \n_f \mathrm{d}s: = \boldsymbol{\mathrm{I}}_p^b ( \u_f \cdot \n_f ),
\label{Ipb}
\end{align}
where $\boldsymbol{\mathrm{I}}_p^b$ denotes a piecewise constant interpolation on $\Gamma$.
By the trace theorem, \eqref{IpProInt} and the last equation,
\begin{align}\label{agamma}
a_\Gamma( \eta_{ph}, [\e])  = & a_{\Gamma}(\eta_{ph}, \,(\boldsymbol{\mathrm{I}}_f \u_f- \u_f) \cdot \n_f
 - (\boldsymbol{\mathrm{I}}_p \u_p -\u_p) \cdot \n_f)
\nn \\
\leq &  C h^{3/2}\|\eta_{ph}\|_{\Gamma}  \|\u_f\|_{H^2}
+g_0 |\langle \eta_{ph}, \, (\boldsymbol{\mathrm{I}}_p \u_p) \cdot \n_f - \u_p \cdot \n_f \rangle |
 \nn\\
\leq & C h^{3/2} \|\eta_{ph}\|_{\Gamma}
\nn \\
\leq & C h^{3/2} (\|\boldsymbol{e}_{ph}\|_{0}  +  \|[\boldsymbol{e}_h]\|_{\Gamma} + h),
\end{align}
where it is noted $\eta_{ph} \in Q_p^h$, and similarly, using \eqref{Ipb}, we have
\begin{align}\label{gamma}
\langle [\e], \, [\e_h] \rangle
&=\langle (\boldsymbol{\mathrm{I}}_f \u_f - \u_f) \cdot \n_f,\, [\e_h] \rangle
-\langle (\boldsymbol{\mathrm{I}}_p \u_p - \u_p) \cdot \n_f,\, [\e_h] \rangle
\nn \\
&\leq |\langle (\boldsymbol{\mathrm{I}}_f \u_f - \u_f) \cdot \n_f,\, [\e_h] \rangle|
+|\langle \boldsymbol{\mathrm{I}}_p^b (\u_f \cdot \n_f) - \u_f \cdot \n_f,\, [\e_h] \rangle|
\leq C h \|[\e_h]\|_{\Gamma}.
\end{align}

To estimate $\|\eta_{fh}\|_0$, take $\boldsymbol{v}_{ph}=0, \, \v_{fh} \in \X_{f,0}^h \subset \X_f^h$ and $q_{ph}=0, \, q_{fh}=0$ in \refe{err-0} again
to obtain
\begin{align}
b_f(\v_{fh}, \, \eta_{fh}) & = a_f(\e_{fh}, \,\v_{fh})
+a_f(\e_f, \, \v_{fh}) -b_f(\v_{fh}, \, \eta_f)
\nn \\
& \le C ( \| \e_{fh} \|_1 + h) \| \v_{fh} \|_1
\nn
\end{align}
which with  the inf-sup condition leads to
\begin{align}\label{etafh}
\| \eta_{fh} \|_0 \le C ( \| \e_{fh} \|_1 + h)  \, .
\end{align}
Similarly, to estimate $\|\eta_{ph}\|_0$, by setting $\boldsymbol{v}_{fh}=0$ and $q_{fh}=0, \, q_{ph}=0$ in \refe{err-0}, we have
\begin{align}
b_p(\v_{ph}, \, \eta_{ph}) - a_\Gamma(\eta_{ph}, \, \v_{ph} \cdot \n_p )
& = a_p(\e_{ph}, \,\v_{ph})
+ \gamma \langle [\e_{h}] , \, \v_{ph}\cdot \n_p \rangle
+ a_p(\e_p, \, \v_{ph})
\nn \\
& -b_p(\v_{ph}, \, \eta_p)
+
a_\Gamma(\eta_p, \, \v_{ph} \cdot \n_p) + \gamma \langle [\e], \, \v_{ph} \cdot \n_p \rangle.
\label{bpinf-sup}
\end{align}
Using the trace theorem and the definition of $\boldsymbol{\mathrm{I}}_p$, we arrive at
\begin{align}
\langle [\e], \, \v_{ph} \cdot \n_p \rangle
&=\langle \,(\boldsymbol{\mathrm{I}}_f \u_f- \u_f) \cdot \n_f
 - (\boldsymbol{\mathrm{I}}_p \u_p -\u_p) \cdot \n_f), \, \v_{ph} \cdot \n_p \rangle
\nn \\
&\leq |\langle (\boldsymbol{\mathrm{I}}_f \u_f - \u_f) \cdot \n_f,\, \v_{ph} \cdot \n_p \rangle|
\nn \\
&\leq C h \|\v_{ph}\|_1. \label{[e]}
\end{align}
By Lemma \ref{3-2}, approximation properties of interpolation/projection operators and \eqref{[e]}, \eqref{bpinf-sup}~can be bounded by
\begin{align}\label{p-p}
\|\eta_{ph}\|_0 + \|\eta_{ph}\|_{\Gamma} \leq
C(\|\boldsymbol{e}_{ph}\|_{0}  +  \|[\boldsymbol{e}_h]\|_{\Gamma} + h).
\end{align}

Combining \eqref{bp}, and \eqref{agamma}--\eqref{etafh}, \eqref{RHS}
follows that
\begin{align}
\mbox{RHS}
\leq Ch(\|\e_{fh}\|_1
+\|\e_{ph}\|_{0}
+ \|[\e_h]\|_{\Gamma} + h) \, .
\nn
\end{align}
Using \eqref{a+gamma}, we have
\begin{align*}
\|\e_{fh}\|_1^2
+\|\e_{ph}\|_{0}^2
+ \|[\e_h]\|_{\Gamma}^2  \leq Ch(\|\e_{fh}\|_1
+\|\e_{ph}\|_{0}
+ \|[\e_h]\|_{\Gamma} + h) \, .
\nn
\end{align*}
Furthermore, by \eqref{etafh} and \eqref{p-p},
the $H^1$-norm estimates are obtained:
\begin{equation}
\begin{split}
&\|\e_{fh}\|_1
+\|\e_{ph}\|_0
+\|[\e_h]\|_{\Gamma}
\leq Ch,\\
&\|\eta_{fh}\|_0 + \|\eta_{ph}\|_0 + \|\eta_{ph}\|_{\Gamma} \leq Ch
\end{split}
\nn
\end{equation}
which with the approximation properties of the corresponding interpolation/projection operators
leads to
\begin{equation}\label{H1-norm}
\begin{split}
&\|\u_f - \wu_{fh}\|_1
+\|\u_p - \wu_{ph}\|_0
+\|[\u - \wu_h]\|_{\Gamma}
\leq Ch,  \\
&\|p_f - \wp_{fh}\|_0 + \|\phi_p - \widetilde{\phi}_{ph}\|_0 + \|\phi_p - \widetilde{\phi}_{ph} \|_{\Gamma} \leq Ch.
\end{split}
\end{equation}

To get the $L^2$-norm estimate for the Stokes velocity, we follow the Aubin--Nitsche duality argument.
Let $\w=(\w_f,\w_p), \r=(r_f,r_p)$, and consider the following dual problem
\begin{align}
-\nabla \cdot \mathbb{T}^{\ast} := -2 \nu \nabla \cdot \mathbb{D}(\w_f)-\nabla r_f=\u_{f} - \wu_{fh}, \label{dualStokes}
~~\nabla \cdot \w_f=0~~~~\text{in}~\Omega_f,\\
\w_p = \boldsymbol{K} \nabla r_p,
~~\nabla \cdot \w_p=0~~~~\text{in}~\Omega_p, \label{dualDarcy}
\end{align}
subject to the interface conditions
\begin{align*}
[\w]&=0~~~~\text{on}~\Gamma,\\
-\n_f \cdot \mathbb{T}^{\ast} \cdot \n_f &= -r_f - 2\nu \n_f \cdot \mathbb{D}(\w_f) \cdot \n_f =-g_0 r_p~~~~\text{on}~\Gamma,\\
-\n_f \cdot \mathbb{T}^{\ast} \cdot \tau_i = -2 \nu \n_f \cdot \mathbb{D}(\w_f) \cdot \tau_i
&= \frac{\alpha \nu \sqrt{d}}{\sqrt{\text{trace}(\Pi)}} \w_f \cdot \tau_i,
~~~~i=1,...,d-1,~~~~\text{on}~\Gamma,
\end{align*}
with homogeneous Dirichlet boundary conditions on the remaining part of the boundary.
The solution of the dual problem enjoys the regularity
\cite{lyu2023regularity}
\begin{align}
\|\w_f\|_2 + \|\w_p\|_1 + \|r_f\|_1 + \|r_p\|_2
\leq \|\u_f - \wu_{fh}\|_0.
\label{reg}
\end{align}
Choosing the test functions $\v_f=\u_f - \wu_{fh}, q_f=p_f-\wp_{fh}, \v_p=\u_p - \wu_{ph}$ and $q_p=\phi_p-\widetilde \phi_{ph}$, we may write
\begin{align*}
\|\u_f-\wu_{fh}\|_0^2
=a(\w, \,\u - \wu_h)
+b(\u - \wu_h, \,\r)
-b(\w, \,\p-\widetilde \p_h)
-a_{\Gamma}(r_p, \,[\u - \wu_h]),
\end{align*}
and by Galerkin orthogonality \refe{Orthogonality},
\begin{align}
\|\u_f-\wu_{fh}\|_0^2
&=a(\w-\boldsymbol{\mathrm{I}}_h \w, \, \u-\wu_h)
+b(\u-\wu_h, \, \r-\mathbb{I}_h \r)
-b(\w-\boldsymbol{\mathrm{I}}_h \w, \, \p-\widetilde \p_h)
\nn \\
&~~~~-a_{\Gamma}(r_p-\mathbb{I}_p r_p,\, [\u-\wu_h])
-a_{\Gamma}(\phi_p -  \widetilde \phi_{ph}, \, [\boldsymbol{\mathrm{I}}_h \w])
-\gamma \langle [\u-\wu_h], \, [\boldsymbol{\mathrm{I}}_h \w]  \rangle,
\nn
\end{align}
where $\boldsymbol{\mathrm{I}}_h \w=(\boldsymbol{\mathrm{I}}_f \w_f, \boldsymbol{\mathrm{I}}_p \w_p)$ and $\mathbb{I}_h \r=(\mathbb{I}_f r_f, \mathbb{I}_p r_p)$.
By noting the approximation properties of these interpolation/projection operators again, the following holds:
\begin{align}
\|\u_f-\wu_{fh}\|_0^2
&\le Ch ( \| \w_f \|_2 \|\u_f-\wu_{fh} \|_1+ \| \w_p \|_1 \|\u_p-\wu_{ph} \|_{0})
\nn \\
&~~~~+ C h \| \u_f - \wu_{fh}\|_1 \|r_f \|_1 +  |b_p(\u_p-\wu_{ph},\, r_p - \mathbb{I}_p r_p)|
\nn \\
&~~~~ + C h \| \w_f \|_2 \| p_f-\widetilde{p}_{fh} \|_0
+ |b_p(\w_p - \boldsymbol{\mathrm{I}}_p \w_p, \, \phi_p- \widetilde \phi_{ph} ) |
+ C h \| r_p \|_{2} \| [\u-\wu_h] \|_\Gamma
\nn \\
&~~~~ -a_{\Gamma}(\phi_p - \widetilde \phi_{ph}, \, [\boldsymbol{\mathrm{I}}_h \w])
-\gamma \langle [\u-\wu_h], \, [\boldsymbol{\mathrm{I}}_h \w]  \rangle
\nn \\
& \le
Ch^2 ( \| \w_f \|_2 + \| \w_p \|_1 +  \|r_f \|_1 +  \| r_p \|_2)
+|b_p(\w_p - \boldsymbol{\mathrm{I}}_p \w_p, \, \phi_p- \widetilde \phi_{ph} ) | \nn\\
&~~~~ + |b_p(\u_p-\wu_{ph},\, r_p - \mathbb{I}_p r_p)|
-a_{\Gamma}(\phi_p - \widetilde \phi_{ph}, \, [\boldsymbol{\mathrm{I}}_h \w])
-\gamma \langle [\u-\wu_h], \, [\boldsymbol{\mathrm{I}}_h \w] \rangle.
\label{l2-1}
\end{align}
Using the divergence-free property of $\w_p$ and the classic commuting diagram property \cite{boffi2013mixed}, we have
\begin{align}\label{b_p}
b_p(\w_p - \boldsymbol{\mathrm{I}}_p \w_p, \, \phi_p- \widetilde \phi_{ph} )
=0,
\end{align}
and
\begin{align}
b_p(\u_p-\wu_{ph},\, r_p - \mathbb{I}_p r_p)&=
b_p(\u_p-\boldsymbol{\mathrm{I}}_p \u_p, \, r_p - \mathbb{I}_p r_p) \nn \\
& \le C h^2 \| \nabla \cdot \u_p \|_1 \| r_p \|_1
\nn \\
& \leq C h^2 \|r_p\|_1.
\end{align}
Then we estimate the last two terms in \refe{l2-1}.
Repeating the procedure from \eqref{agamma} to \eqref{gamma} and noting that $[\boldsymbol{\mathrm{I}}_h \w]=[\boldsymbol{\mathrm{I}}_h \w - \w]$ on $\Gamma$,  the second last term is bounded by
\begin{align}
a_\Gamma(\phi_p - \widetilde \phi_{ph}, [\boldsymbol{\mathrm{I}}_h \w] )
& = a_\Gamma(\phi_p - \widetilde \phi_{ph},  \, (\boldsymbol{\mathrm{I}}_f \w_f- \w_f) \cdot \n_f)
- a_{\Gamma}(\phi_p - \widetilde \phi_{ph},  \,  (\boldsymbol{\mathrm{I}}_p \w_p  - \w_p) \cdot \n_f )
\nn \\
& \le C h \| \phi_p - \widetilde \phi_{ph} \|_\Gamma  \| \w_f \|_{2}
+g_0 | \langle \phi_p - \mathbb{I}_p \phi_p,  \,  \boldsymbol{\mathrm{I}}_p^b(\w_f\cdot \n_f) - (\w_f\cdot \n_f)  \rangle |
\nn \\
& \le Ch \| \phi_p - \widetilde \phi_{ph} \|_{\Gamma} \| \w_f \|_2
+ Ch^2 \| \w_f \|_{2}
\nn \\
& \le Ch^2  \| \u_f - \wu_{fh} \|_{L^2(\Omega_f)},
\label{l2-2}
\end{align}
where we have used the trace theorem, \refe{Ipb} and  \eqref{H1-norm} and noted the regularity bound \refe{reg}.
For the last term in \refe{l2-1}, similarly,
\begin{align}
\langle [\u-\wu_h], \,[\boldsymbol{\mathrm{I}}_h \w]  \rangle
& =
\langle [\u -\wu_{h}]\cdot \n_f, \, (\boldsymbol{\mathrm{I}}_f \w_f - \w_f) \cdot \n_f \rangle
+ \langle [\u -\wu_{h}]\cdot \n_f,  \, (\w_p - \boldsymbol{\mathrm{I}}_p \w_p)\cdot \n_f \rangle
\nn \\
& \le Ch  \| [\u -\wu_{h}]  \|_{\Gamma} \| \w_f \|_2
\nn \\
& \le C h^2 \| \w_f \|_2,
\label{l2-3}
\end{align}
where \eqref{H1-norm} has been used again.
By the estimates \refe{b_p}-\refe{l2-3}, \refe{l2-1} reduces to
\begin{align}
\|\u_f-\wu_{fh}\|_0^2
&\leq Ch^2
(\|\w_f\|_2 + \|\w_p\|_1+\|r_f\|_1 + \|r_p\|_{2}
+ \|\u_f - \wu_{fh}\|_0 )
\nn\\
&\leq C h^2 \|\u_f - \wu_{fh}\|_0,
\end{align}
using the regularity bound \refe{reg} and the $H^1$-norm estimates \eqref{H1-norm}.
The results \refe{err-sd-1}-\refe{err-sd-2} follow immediately.

Next, we establish the $H^{-1}$-norm estimate.
For any given $\psi_p$, introduce the following dual problem
\begin{align}
-\nabla \cdot \mathbb{T}^{\ast} := -2 \nu \nabla \cdot \mathbb{D}(\w_f)-\nabla r_f=0, \label{dualStokes2}
~~\nabla \cdot \w_f=0~~~~\text{in}~\Omega_f,\\
\w_p = \boldsymbol{K} \nabla r_p,
~~-\nabla \cdot \w_p=\psi_p~~~~\text{in}~\Omega_p, \label{dualDarcy2}
\end{align}
subject to the interface conditions
\begin{align*}
[\w]&=0~~~~\text{on}~\Gamma,\\
-\n_f \cdot \mathbb{T}^{\ast} \cdot \n_f &= -r_f - 2\nu \n_f \cdot \mathbb{D}(\w_f) \cdot \n_f =-g_0 r_p~~~~\text{on}~\Gamma,\\
-\n_f \cdot \mathbb{T}^{\ast} \cdot \tau_i = -2 \nu \n_f \cdot \mathbb{D}(\w_f) \cdot \tau_i
&= \frac{\alpha \nu \sqrt{d}}{\sqrt{\text{trace}(\Pi)}} \w_f \cdot \tau_i,
~~~~i=1,...,d-1,~~~~\text{on}~\Gamma,
\end{align*}
with homogeneous Dirichlet boundary conditions on the rest part of the boundary. Again, the solution of the dual problem enjoys the regularity
\begin{align}
\|\w_f\|_2 + \|\w_p\|_1 + \|r_f\|_1 + \|r_p\|_2
\leq \|\psi_p\|_{0}.
\label{reg2}
\end{align}
Choosing the test functions $\v_f=\u_f - \wu_{fh}, q_f=p_f-\wp_{fh}, \v_p=\u_p - \wu_{ph}$ and $q_p=\phi_p-\widetilde \phi_{ph}$, we may write
\begin{align*}
(\phi_p- \widetilde \phi_{ph},\psi_p)
=a(\w, \, \u-\wu_h)
+b(\u-\wu_h, \, \r)
-b(\w, \, \p - \widetilde \p_{h})
- a_{\Gamma}(r_p,\, [\u-\wu_h]),
\end{align*}
and by noting Galerkin orthogonality,
\begin{align}
(\phi_p- \widetilde \phi_{ph},\, \psi_p)
&=a(\w -\boldsymbol{\mathrm{I}}_h \w, \, \u -\wu_{h})
+b(\u-\wu_h, \, \r-\mathbb{I}_h \r)
-b(\w-\boldsymbol{\mathrm{I}}_h \w, \, \p- \widetilde \p_{h}) \nn\\
&~~~~-a_{\Gamma}(r_p - \mathbb{I}_p r_p, \, [\u-\wu_h])
-a_{\Gamma}(\phi_p - \widetilde \phi_{ph}, \, [\boldsymbol{\mathrm{I}}_h \w])
-\gamma \langle [\u-\wu_h], \, [\boldsymbol{\mathrm{I}}_h \w]  \rangle . \label{H-1}
\end{align}
Using similar estimations as in \eqref{l2-1},
with the regularity in \eqref{reg2}, we obtain
\begin{align}
(\phi_p- \widetilde \phi_{ph},\, \psi_p) \le C h^2 \| \psi_p \|_1,
\nn
\end{align}
which leads to the $H^{-1}$-norm estimate
\begin{align*}
\|\phi_p - \widetilde \phi_{ph}\|_{H^{-1}}
\leq \sup_{0 \neq \psi_p \in H^1_{0,D}} \frac{(\phi_p - \widetilde \phi_{ph}, \,  \psi_p)}{\|\psi_p\|_1}
\leq C h^2.
\end{align*}
Moreover, taking the difference quotient of the equation \refe{err-0} provides
\begin{align}
&a(d_t \e_{h}^{n+1}, \,\v_h)
-b(\v_h, \, d_t \eta_h^{n+1})
+b(d_t \e_h^{n+1}, \, \q_h)
+a_{\Gamma}(d_t \eta_{ph}^{n+1}, \,  [\v_h])
-a_{\Gamma}(q_{ph},  \, [d_t \e_h^{n+1}])
\nn \\
&+\gamma \langle [d_t \e_h^{n+1}], \,  [\v_h] \rangle
\nn \\
&=-\big( a(d_t \e^{n+1}, \, \v_h)
-b(\v_h, \, d_t \eta^{n+1})
+b(d_t \e^{n+1}, \, \q_h)
+a_{\Gamma}(d_t \eta_{p}^{n+1}, \, [\v_h])
\nn \\
&\quad
-a_{\Gamma}(q_{ph}, \, [d_t \e^{n+1}])
+\gamma \langle [d_t \e^{n+1}], \, [\v_h] \rangle   \big).
\label{err-t}
\end{align}
Using the same approach as above, we can prove the inequality \refe{err-sd-3}.
The proof of Theorem \ref{T3-1} is completed. \endproof

\section{Numerical results}
In this section, we present our numerical results by two examples.
The first one is an artificial model with a smooth exact solution to confirm our theoretical analysis, and the second is a practical model arising from the standard technique extracting oil/gas from a reservoir.
All computations were performed using the FreeFem++ software package \cite{Hechet2010}.
\vskip0.1in

{\bf Example 5.1.} We consider an artificial example defined in \refe{e-f}--\refe{initial}, where
the computational domain~$\Omega$~is composed of~$\Omega_f=(0,1)\times (0,1)$~and~$\Omega_p=(0,1)\times (1,2)$, with the interface~$\Gamma=(0,1)\times \{1\}$ and $\Gamma_{p}^D=(0,1) \times \{0\}$.
Setting the hydraulic parameters $k_1 = k_2=1$ and $\nu=1, S_0=1, g_0=1, \alpha=1$ for simplicity, the analytical solutions satisfying the time-dependent Stokes--Darcy model are as follows:
\begin{align*}
&\boldsymbol{u}_f =\Big[ \big( x^2(y-1)^2+y \big)\cos(t),~ -\frac{2}{3}x(y-1)^3\cos(t)+\big( 2-\pi \sin(\pi x) \big) \cos(t) \Big ]^{T},\\
&\boldsymbol{u}_p=\Big[ \pi^2 \cos(\pi x)(1-y-\cos(\pi y))\cos(t),~
(2-\pi \sin(\pi x))(1-\pi \sin(\pi y)) \cos(t) \Big]^{T},\\
&p_f = \big( 2-\pi \sin(\pi x) \big)  \sin(0.5 \pi y) \cos(t),\\
&\phi_p = \big( 2-\pi \sin (\pi x) \big) \big( 1-y-\cos(\pi y) \big) \cos(t).
\end{align*}
The initial conditions, boundary conditions and the forcing terms are determined by the analytical solutions.

We first solve the coupled Stokes--Darcy problem by the proposed algorithm with
the lowest-order fully-mixed FE approximations (MINI/RT0-DG0).
To show the convergence rates, we choose $\tau=h^2, \gamma=1.0$ in our numerical simulations.
The numerical results at $t=1.0$ with $h=1/4, 1/8, 1/16, 1/32, 1/64$ are presented in Table \ref{T1}.
These results clearly demonstrate that the algorithm achieves second-order accuracy in the $L^2$-norm for the velocity of Stokes flows, and first-order accuracy for the velocity/pressure of Darcy flow in the spatial direction.
For the coupled model, the one-order lower approximation to $(\phi_p, \u_p)$ in Darcy region
does not affect the accuracy of the velocity with a higher finite element approximation in Stokes region, which are optimal for all physical variables and consistent with our theoretical results.

\begin{table}[htbp!]
\caption{\label{T1}$L^2$ errors of the decoupled scheme with non-uniform $\mbox{P1b-P1/RT0-DG0}$ approximations}
\centering
\begin{tabular}{ccccccc}
\hline
$h$ &  $\|\u_f(t_N) -\u_{fh}^N\|_0$ &
Rate &
$\|\u_p(t_N) -\u_{ph}^N\|_0$&
Rate &
$\|\phi_p(t_N) - \phi_{ph}^N\|_0$&
Rate \\
\hline
$\frac{1}{4}$  & 6.693E-02    & --      & 2.975E-01   & --     & 4.432E-01    &--  \\
$\frac{1}{8}$  & 1.691E-02    & 1.98    & 1.489E-01   & 1.00   & 2.218E-01    &1.00 \\
$\frac{1}{16}$ & 4.237E-03    & 2.00    & 7.445E-02   & 1.00   & 1.109E-01    &1.00 \\
$\frac{1}{32}$ & 1.060E-03    & 2.00    & 3.722E-02   & 1.00   & 5.543E-02    &1.00 \\
$\frac{1}{64}$ & 2.649E-04    & 2.00    & 1.861E-02   & 1.00   & 2.771E-02    &1.00   \\
\hline
\end{tabular}
\end{table}

\begin{table}[htbp!]
\caption{\label{T2}$L^2$ errors of the decoupled scheme with non-uniform $\mbox{P2-P1/RT0-DG0}$ approximations}
\centering
\begin{tabular}{ccccccc}
\hline
$h$ &  $\|\u_f(t_N) -\u_{fh}^N\|_0$ &
Rate &
$\|\u_p(t_N) -\u_{ph}^N\|_0$&
Rate &
$\|\phi_p(t_N) - \phi_{ph}^N\|_0$&
Rate \\
\hline
$\frac{1}{4}$  & 3.211E-03    & --      & 2.985E-01   & --     & 4.774E-01    &--  \\
$\frac{1}{8}$  & 4.668E-04    & 2.78    & 1.490E-01   & 1.00   & 2.269E-01    &1.07 \\
$\frac{1}{16}$ & 1.037E-04    & 2.17    & 7.446E-02   & 1.00   & 1.115E-01    &1.02 \\
$\frac{1}{32}$ & 2.680E-05    & 1.95    & 3.722E-02   & 1.00   & 5.551E-02    &1.01 \\
$\frac{1}{64}$ & 6.942E-06    & 1.95    & 1.861E-02   & 1.00   & 2.772E-02    &1.00   \\
\hline
\end{tabular}
\end{table}

Moreover, we solve the coupled model using the Taylor--Hood and Raviart--Thomas finite element pair $\mbox{P2-P1/RT0-DG0}$, which is also non-uniform approximations.
We present our numerical results for $\tau=h^3$ in Table \ref{T2} and observe that the $L^2$-norm error for the Stokes velocity is of order $O(h^2)$ in spatial direction, which is one-order lower than that for interpolation, but optimal for the coupled model.
It
shows that in this case, the lower-order approximation to the Darcy flow does pollute the accuracy of the numerical velocity for Stokes flow.


\vskip0.1in

{\bf Example 5.2. (Horizontal open-hole completion with vertical wellbore)}
We consider an example in petroleum engineering:
a horizontal open-hole completion with a vertical wellbore, which
is a standard technique for extracting oil/gas from a reservoir \cite{al2007innovative, wan2011advanced}.
To enhance reservoir recovery, water, superheated steam, or specialized chemicals (such as superhydrophobic proppants used in shale oil/gas extraction to reduce water blockage and improve hydrocarbon flow paths) are often injected through vertical wellbores to maintain pressure and improve the flow of hydrocarbons to the production well
\cite{donaldson2014hydraulic}.
To optimize oil recovery techniques and improve production efficiency,
numerous numerical simulations are performed to thoroughly assess the reservoir's production potential.
The fluid velocities within both the wellbore and the reservoir are critical physical components that significantly influence recovery efficiency \cite{tiab2024petrophysics}.

However, the reservoir area is extremely large compared to the pipeline's cross-sectional area, and the extended production period makes it necessary to efficiently simulate the recovery process while saving computational costs.
Moreover, the irregularity of the actual reservoir media makes lower-order elements more suitable for simulation.
Therefore, the lowest-order fully-mixed FE scheme can greatly improve the efficiency of simulation.
Moreover, the non-pollution of the lowest-order elements ensures optimal accuracy in the fluid flow velocity within the production wellbore.

\begin{figure}[htbp!]
\begin{centering}
\begin{subfigure}[t]{0.3\textwidth}
\centering
\includegraphics[width=1.55\textwidth]{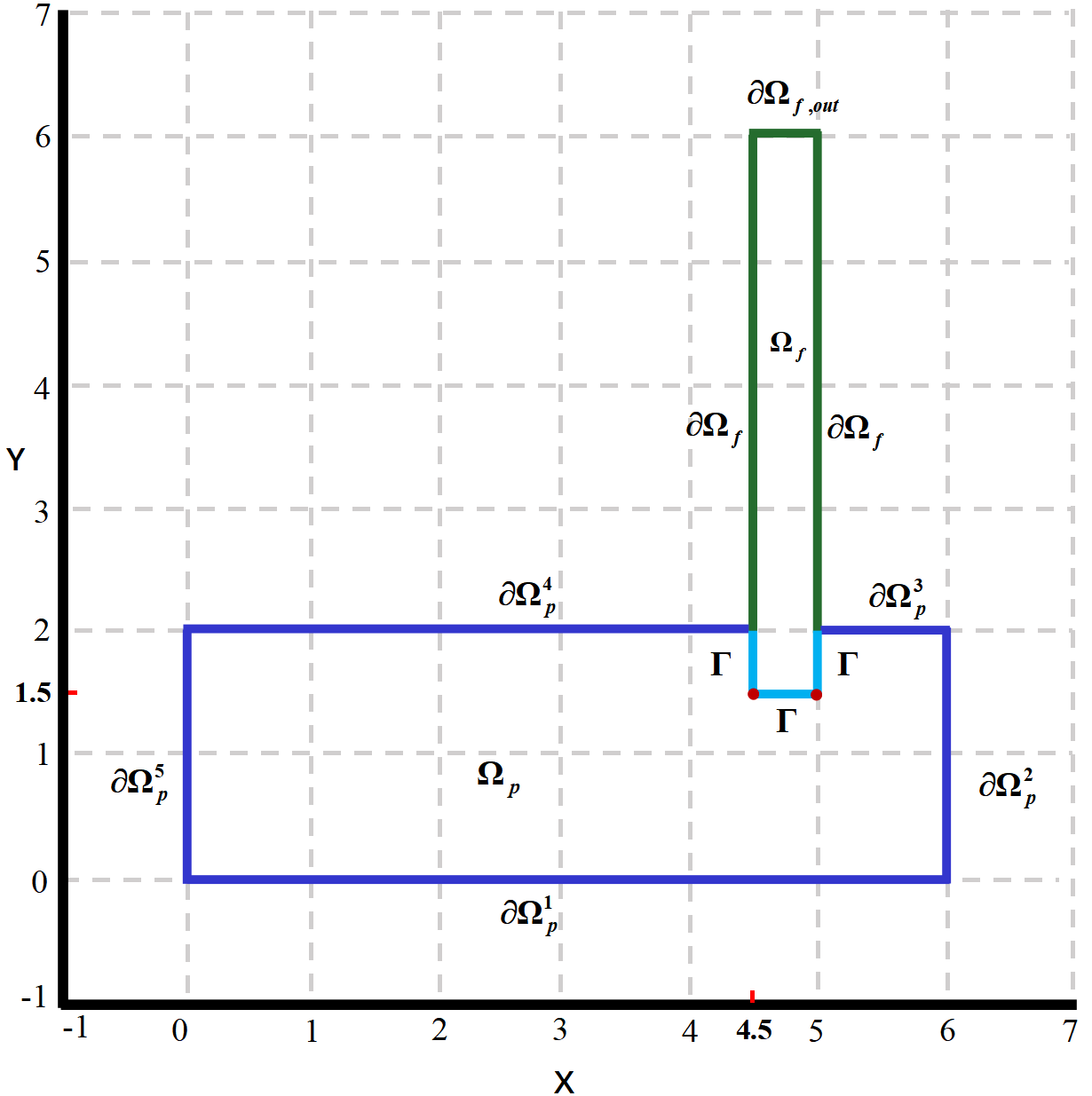}
\end{subfigure}
\hspace{30mm}
\begin{subfigure}[t]{0.33\textwidth}
\centering
\includegraphics[width=1.5\textwidth]{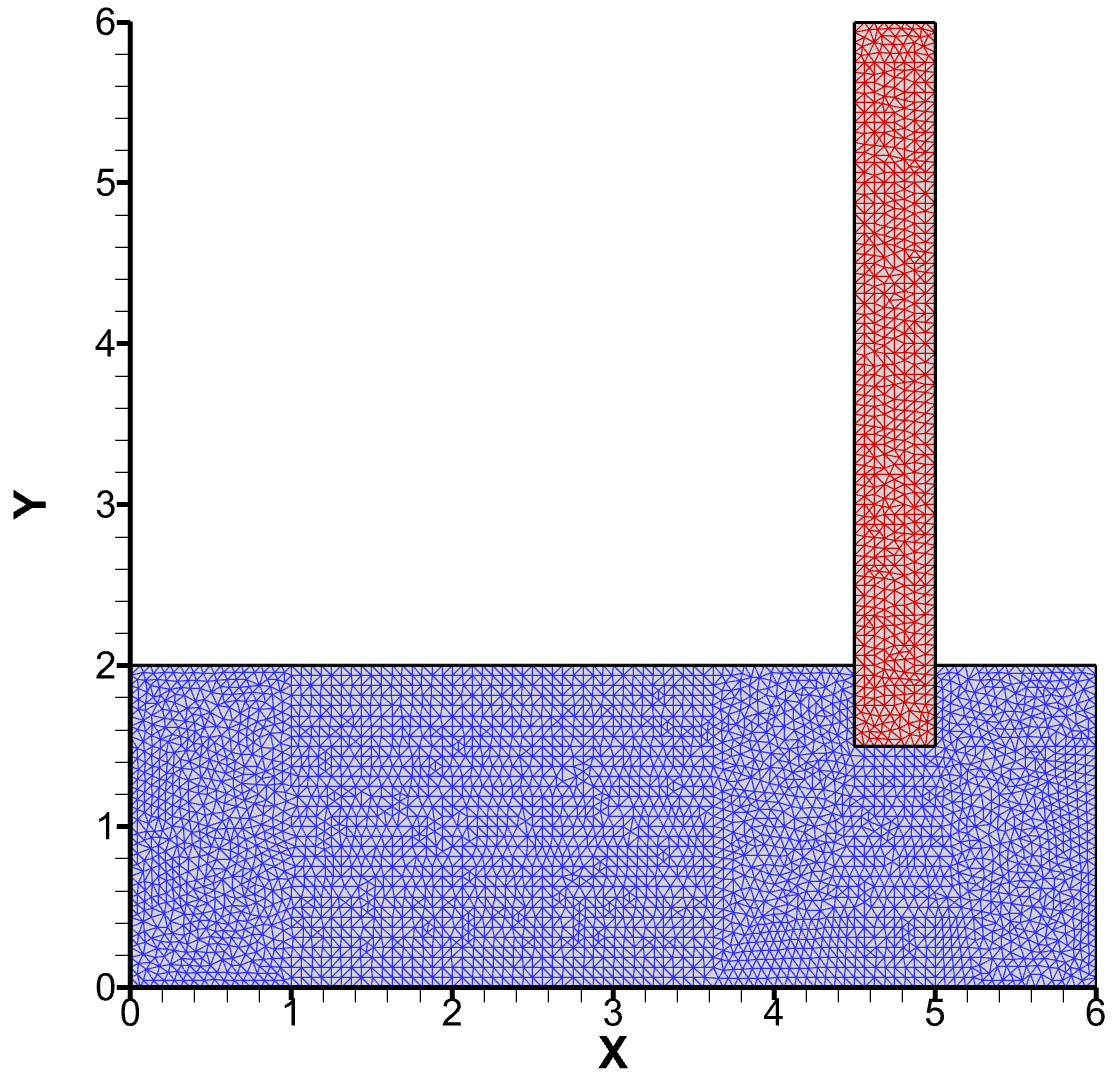}
\end{subfigure}
\end{centering}
\caption{\small{\label{Fig3} Left: the geometric location of the simulation domain; Right: quasi-uniform grid.}}
\end{figure}

As shown in Figure \ref{Fig3}, the geometric location of the simulation is depicted with the coordinate axes, and a quasi-uniform grid is used.
The vertical pipeline functions for both injection and production, but this simulation focuses on the oil/gas production process.
The Neumann boundary condition $(-p_f \mathbb{I} + \nu \nabla \u_f) \cdot \n_f = 0$ is imposed on the top of vertical wellbore.
On the remaining boundaries of the vertical wellbore,
$\partial \Omega_f = \{(x,y): x=4.5, 2 \leq y \leq 6\} \cup \{(x,y): x=5, 2 \leq y \leq 6\}$,
the non-slip condition $\boldsymbol{u}_f=(0,0)$ is imposed.
The porosity flow velocity $\u_p$ on the boundaries $\partial \Omega_p^i(i=1,2,3,4,5)$ are assumed as $(0,\Theta_1), (-\Theta_1,0), (0,-\Theta_1), (0,-\Theta_1), (\Theta_1,0), \Theta_1=1.0$, respectively.

\begin{figure}[htbp!]
\begin{centering}
\begin{subfigure}[t]{0.3\textwidth}
\centering
\includegraphics[width=1.65\textwidth]{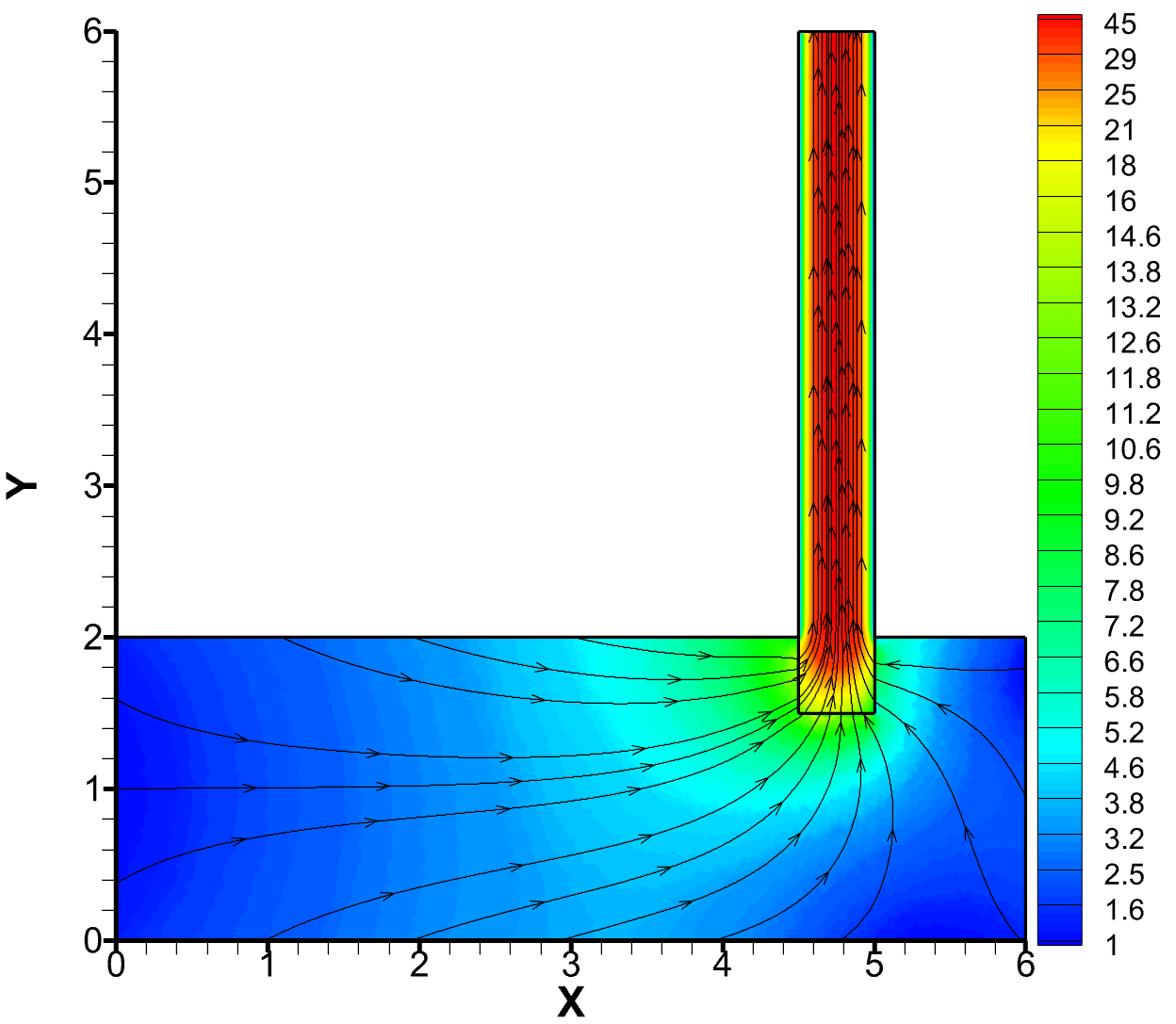}
\end{subfigure}
\hspace{30mm}
\begin{subfigure}[t]{0.3\textwidth}
\centering
\includegraphics[width=1.65\textwidth]{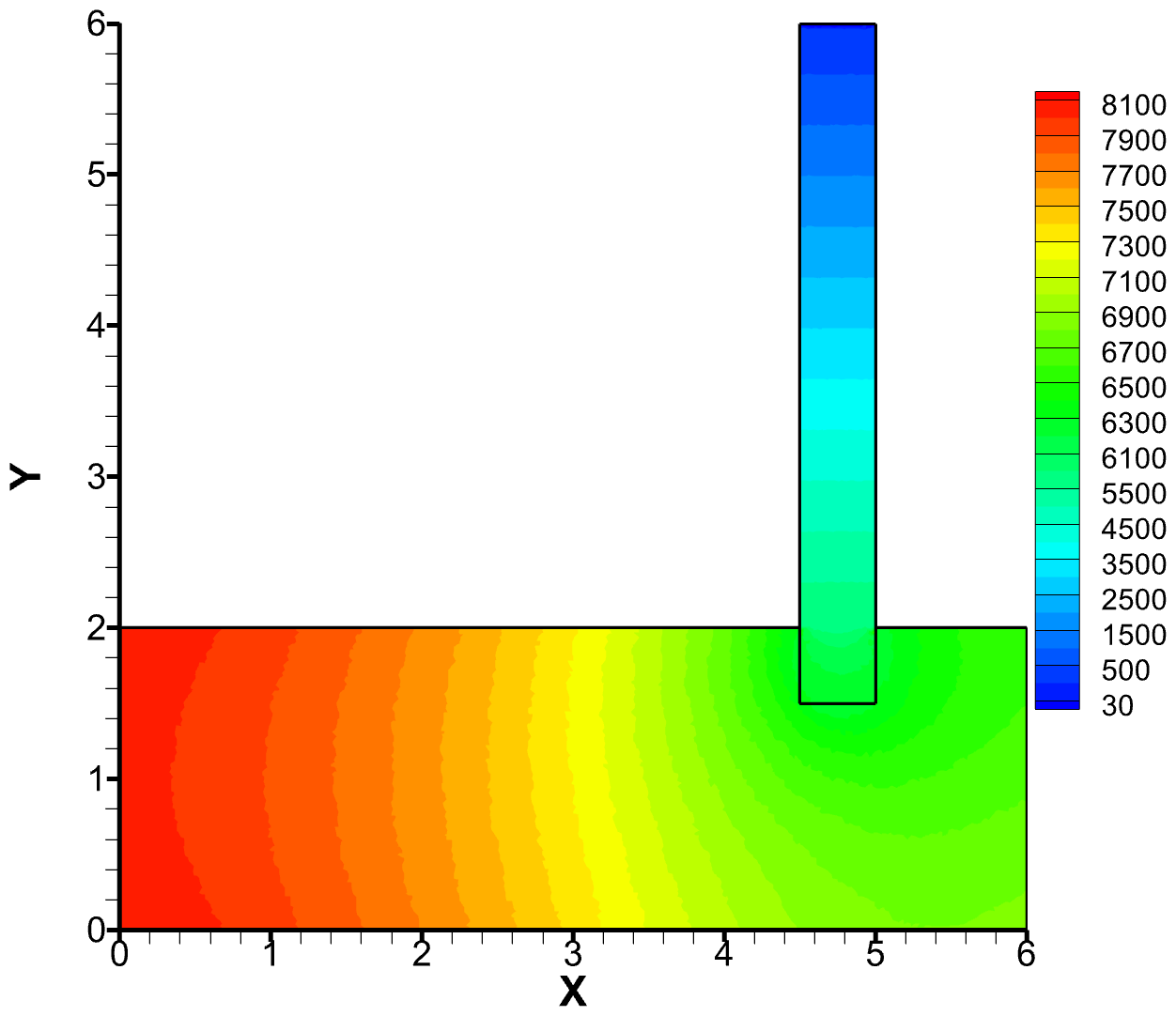}
\end{subfigure}
\end{centering}
\caption{\label{Fig4}The flow speed, streamlines and pressure around horizontal open-hole attached with vertical wellbore completion for reservoir velocity $\Theta_1=1$. Left: flow speed streamlines; Right: pressure.}
\end{figure}

The model parameters are chosen as in \cite{mahbub2023uncoupling},
\begin{align}
& \nu=1.0, \quad S_0=10^{-4}, \quad g_0=1, \quad \alpha=1,
\nn \\
& k_1=k_2=10^{-2},  \quad \boldsymbol{f}_f=0, \quad f_p=0
\nn
\end{align}
with a penalty parameter $\gamma=10$.
The mesh size and time step size are set to $h=1/32$ and $\tau=0.001$, respectively.
We also use $\mbox{P1b-P1/RT0-DG0}$ to approximate the velocity and pressure in the wellbore region and reservoir region.
The numerical results at $t=3.0$ are presented in Figure \ref{Fig4}, where warmer colors indicate higher flow speeds and pressures, and the lines with arrows represent streamlines.
The velocity in the Darcy region shows the flow distribution and the simulation of fluid flow trajectories from the reservoir to the wellbore, which helps optimize the placement of wells.
As the primary storage medium for oil and gas, the porous block exhibits higher pressure, as expected, driving fluid into the wellbore where pressure is comparatively lower.
This pressure difference increases the flow speed in the vertical production pipeline.

\begin{figure}[htbp!]
\begin{centering}
\begin{subfigure}[t]{0.3\textwidth}
\centering
\includegraphics[width=1.65\textwidth]{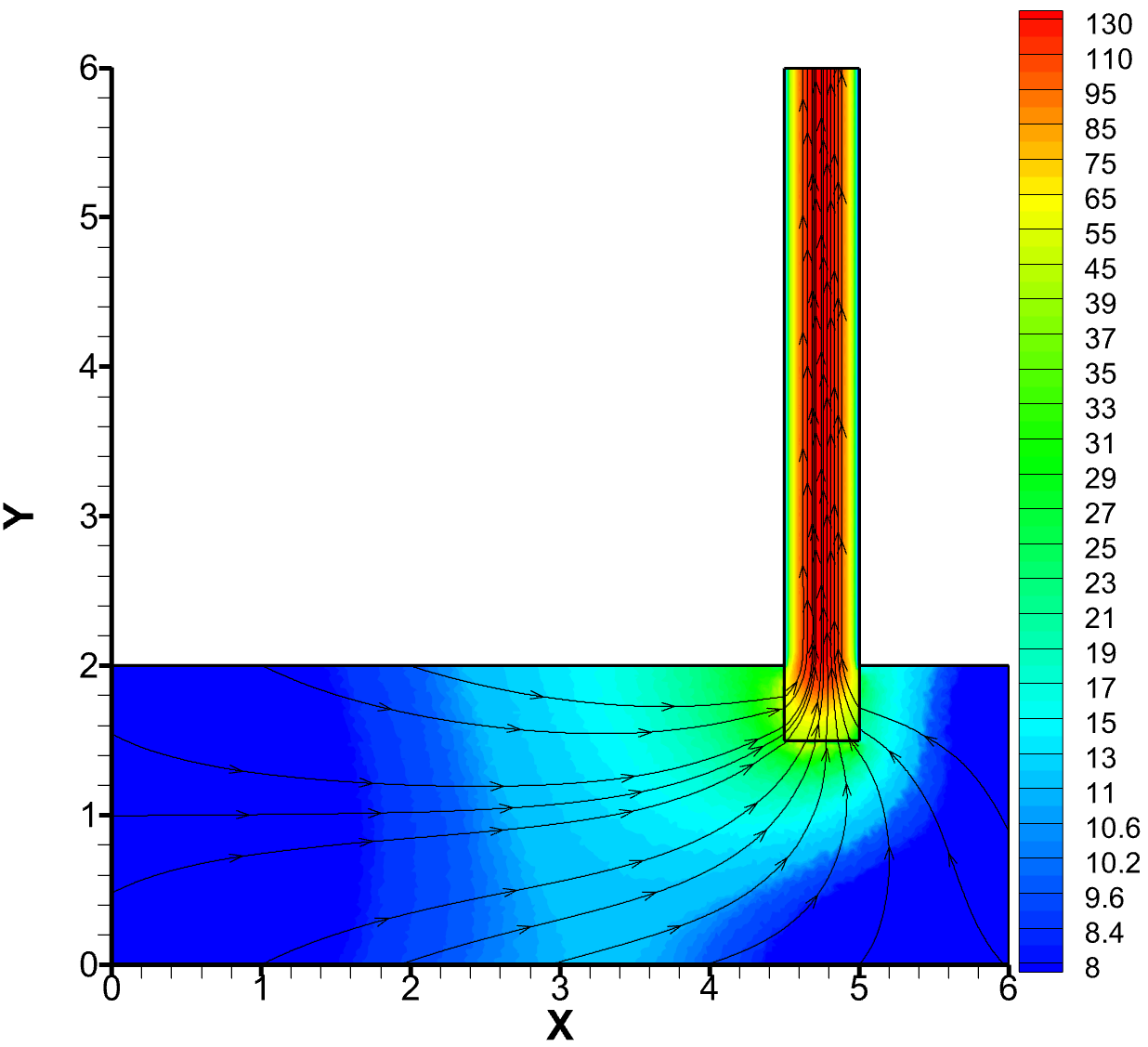}
\end{subfigure}
\hspace{30mm}
\begin{subfigure}[t]{0.3\textwidth}
\centering
\includegraphics[width=1.65\textwidth]{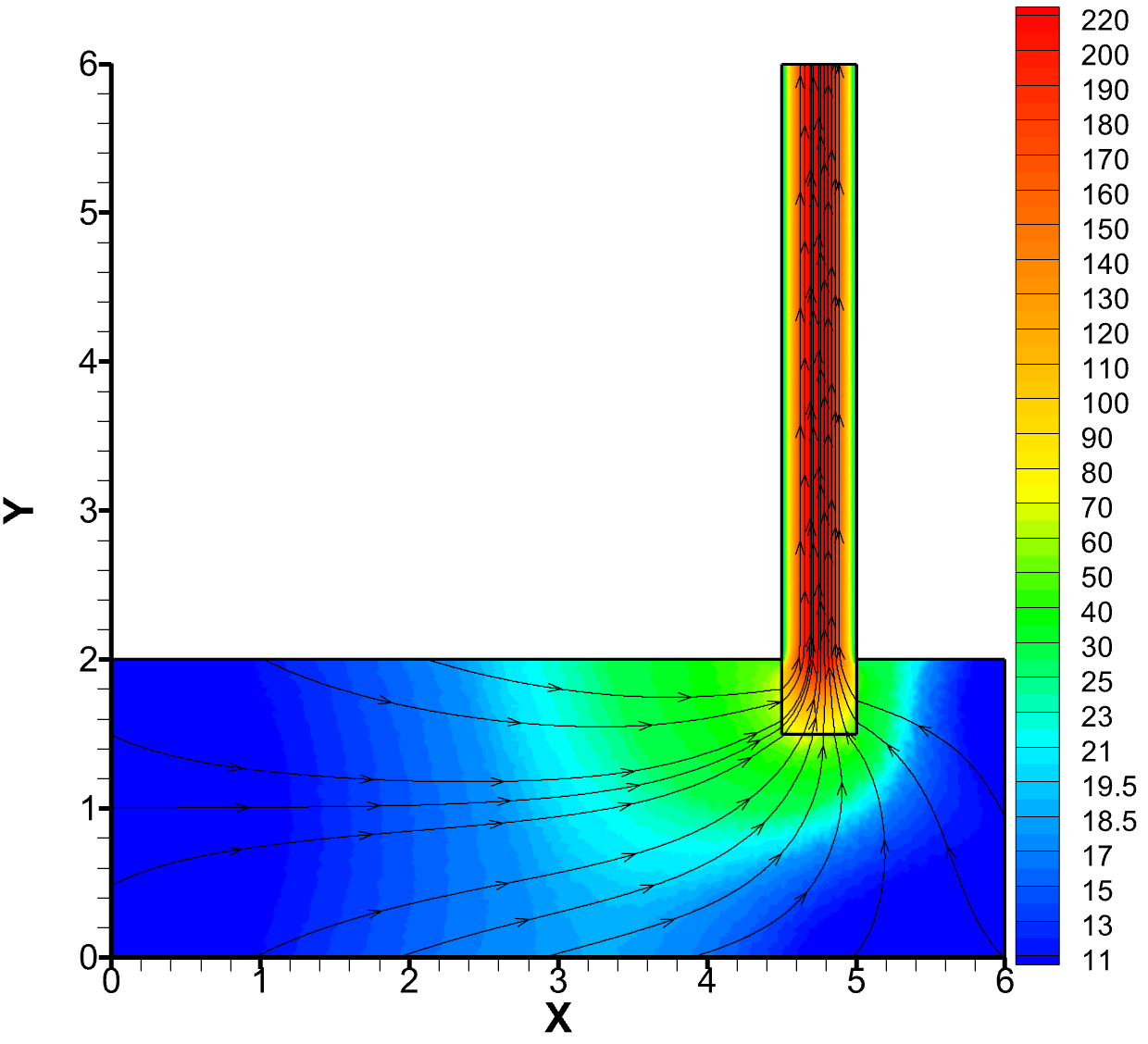}
\end{subfigure}
\end{centering}
\caption{\label{Fig5}The flow speed and streamlines around horizontal open-hole attached with vertical production wellbore completion for different reservoir velocities. Left: $\Theta_2=3$; Right: $\Theta_3=5$.}
\end{figure}

As a comparative numerical experiment, we adjust the reservoir velocity $\u_p$ from $\Theta_2=3.0$ to $\Theta_3=5.0$ along the boundaries $\partial \Omega_p^i$,
with these variations influenced by factors such as the chemical materials mentioned above.
Observing the maximum velocity increase from 130 to 220 in the vertical wellbore in Figure \ref{Fig5}, it can be seen that as the reservoir velocity increases, the recovery rate within the wellbore also improves.
In this way, one can explore the impact of different factors on fluid velocity in the reservoir region and further investigate their effects on the recovery rate.

\section{Declarations}
\textbf{Conflict of interest}
The authors have no Conflict of interest to declare that are relevant to the content of this article.

\section*{Conclusion}
In this paper, we have proposed a decoupled fully-mixed FEM for solving dynamic coupled Stokes--Darcy model.
By employing the lowest-order combined non-uniform approximations (MINI/RT0-DG0),
our approach significantly improves computational efficiency, making it particularly suitable for industrial applications.
Moreover, we present a rigorous error estimate for the proposed algorithm, which is optimal in traditional sense for each physical components involved although the combined approximation is non-uniform. This demonstrates that the lower-order approximation to the Darcy flow does not pollute the accuracy of the numerical velocity for the Stokes flow.
The method and analysis presented in this paper can be extended to other coupled Stokes--Darcy models \cite{badia2009coupling, chen2017uniquely,
vassilev2009coupling}, as well as more general coupled transport equations, with either penalty or Lagrange multiplier approach.

\bibliographystyle{abbrv}
\bibliography{bibfile}

\end{document}